\documentclass[11pt,twoside,a4paper]{article}
\usepackage{mathrsfs}
\usepackage{mathrsfs}
\usepackage{mathrsfs}
\usepackage{mathrsfs}
\usepackage{mathrsfs}
\usepackage{mathrsfs}
\usepackage{mathrsfs}

\usepackage{amsfonts}
\usepackage{amsmath,amssymb}
\usepackage{mathrsfs}
\usepackage{hyperref}

\newtheorem{thm}{Theorem}[section]

\newtheorem{lem}[thm]{Lemma}
\newtheorem{prop}[thm]{Proposition}

\newtheorem{rem}[thm]{Remark}

\numberwithin{equation}{section}\allowdisplaybreaks

\def\le{\leqslant}
\def\ge{\geqslant}
\def\leq{\leqslant}
\def\geq{\geqslant}

\begin{document}

\title{ {\bf \Large Inviscid limit for the derivative
 Ginzburg--Landau equation with small data in
higher spatial dimensions }}
\author{\textsc{\normalsize Lijia Han\footnote{Email:
 hljmath@gmail.com}} \quad\textsc{\normalsize Baoxiang Wang\footnote{
 Email: wbx@math.pku.edu.cn}}\quad \textsc{\normalsize Boling Guo}\footnote{gbl@mail.iapcm.ac.cn}\\
 $^{*,  \ddag}$\textit{\footnotesize {Institute of
Applied Physics and Computational Mathematics,  PO Box 8009, Beijing
100088, China}}\\
[1ex] $^{\dag}$\textit{\footnotesize {LMAM, School of Mathematical Sciences, Peking University,
Beijing 100871, China.}} }
\date{}
\maketitle

\begin{minipage}{13.5cm}
\footnotesize \bf Abstract. \rm  \quad We study the Cauchy problem
for derivative Ginzburg--Landau equation $u_t=  (\nu + i)\triangle u
+\overrightarrow{\lambda_1}\cdot\nabla(|u|^2u) +
(\overrightarrow{\lambda_2}\cdot\nabla u)|u|^2+ \alpha
|u|^{2\delta}u$, where $\delta\in \mathbb{N}$,
$\overrightarrow{\lambda_1}, \overrightarrow{\lambda_2}$ are complex
constant vectors, $ \nu \in [0, 1]$, $\alpha \in \mathbb{C}$. For $
n\geq 3$, we show that it is uniformly global wellposed for all $\nu
\in [0,1]$ if initial data $u_0$ belong to modulation space $M_{2,
1}^{s}$ ($s>3$) and $\|u_0\|_{L^2} \ll 1$. Moreover, we show that
its solution will converge to that of the derivative Schr\"odinger
equation in $C(0,T; L^2)$ if $\nu\rightarrow 0$ and $u_0 \in M_{2,
1}^{4}$. For $n= 2$, we obtain the local well-posedness results and
inviscid limit with  the Cauchy data in $M_{1, 1}^{s}$ ($s>3$) and
$\|u_0\|_{L^1} \ll 1$.

\vspace{10pt}

\bf 2000 Mathematics Subject Classifications. \rm  35 Q 55, 35 K 55.\\

\bf Key words and phrases. \rm Derivative Ginzburg-Landau equation;
Derivative Schr\"odinger equation; Inviscid limit, Modulation spaces.

\end{minipage}

\section{introduction}

In this paper, we consider the Cauchy problem for the derivative
complex Ginzburg--Landau (DCGL) equation:
\begin{align}
u_t=  (\nu + i)\triangle u
+\overrightarrow{\lambda_1}\cdot\nabla(|u|^2u) +
(\overrightarrow{\lambda_2}\cdot\nabla u)|u|^2+ \alpha |u|^{2\delta}u,  \quad u(0,x)=u_0(x), \label{GL}
\end{align}
where $u$ is a complex valued function of $(t, x)\in
\mathbb{R}^+\times \mathbb{R}^n$, $\mathbb{R}^+ = [0, +\infty]$; $\nu>0$,  $\alpha \in \mathbb{C}$,
$\delta\in \mathbb{N}$, $\overrightarrow{\lambda_1}$ and
$\overrightarrow{\lambda_2}$ are complex vectors.

The DCGL equation \eqref{GL} arises as the envelope equation for a
weakly subcritical bifurcation to counter-propagating waves, and it
is also important for a number of physical systems including the
onset of oscillatory convection in binary fluid mixture; cf.
\cite{2}.  In the case of one or two dimensions, the global
existence of solutions, finite dimensional global attractors, Gevery
regularity of solutions have been studied extensively
 for equation \eqref{GL}; cf. \cite{JDE, 12, GW1, LG, WGZ}.  Taking $\nu=0$, \eqref{GL} can be written as
\begin{align}
u_t - i \triangle u = \overrightarrow{\lambda_1}\cdot\nabla(|u|^2u)
+
(\overrightarrow{\lambda_2}\cdot\nabla u)|u|^2+ \alpha |u|^{2\delta}u, \quad u(0,x)=u_0(x),\label{nls}
\end{align}
which is the well-known derivative nonlinear Schr\"odinger equation
(DNLS). There are some recent works which have been devoted to
 equation \eqref{nls}; cf. \cite{14, 15, 18,
19, 22,  WH}. N. Hayashi and Ozawa in \cite{NH} proposed the method
of gauge transformation which is useful to avoid the loss of
derivatives for equation \eqref{nls} in one spacial dimension.

A natural question between Eqs. \eqref{GL} and  \eqref{nls} is the
inviscid limit. Let $u$ and $v$ be the solutions of the
Cauchy problems of Eqs. \eqref{GL} and
\eqref{nls}, respectively. Does $u$ converge to $v$ as the parameter
 $\nu $ tends to $0$?

When $
\overrightarrow{\lambda_1}=\overrightarrow{\lambda_2}=0$,
Eq. \eqref{GL} can be rewritten as
\begin{align}
u_t=  (\nu + i)\triangle u + \alpha |u|^{2\delta}u,\quad\quad
u(0,x)=u_0(x),\label{GL0}
\end{align}
which is the well-known complex Ginzburg--Landau equation. Eq.
\eqref{GL0} is an important model equation in the description of
spatial pattern formation and of the onset of instabilities in
nonequilibrium fluid dynamical systems; cf. \cite{4}. For Eq.
\eqref{GL0}, there are some recent results devoted to the  global
well-posedness and limit behavior, see Ginibre and Velo \cite{GV2},
Wu \cite{WU}, Bechouche and Jungel \cite{BJ}, Wang \cite{Wang},
Machihara and Nakamura \cite{MN}, Wang and Huang \cite{HH}.

For the derivative complex Ginzburg--Landau equation \eqref{GL},
using Bourgain's $X^{s,b}$ method,  Huo and Jia \cite{Huo}  obtained
the inviscid limit for the solutions in $C([0, T]; H^s)$  ($s>1/2$)
in one spatial dimension, where the bilinear estimate condition
 $2\vec{\lambda}_1 + \vec{\lambda}_2 =0$ and some energy estimate conditions on coefficients
 and $\|u_0\|_{L^2} \ll 1$ are required. B. Wang and Y. Wang in \cite{Wa2} also considered
 the inviscid limit for the solutions, when initial data belong to $\dot{H}^3\cup \dot{H}^{-\frac{1}{2}}$, in one spacial dimension.
 As far as the authors can see, there are no
result on the inviscid limit of Eq. \eqref{GL} in high dimension
case $n\geq 2$.

It was well known that $H^{s+\epsilon +n/2}\subset M_{2, 1}^s
\subset H^s$, for $\forall \epsilon>0$. In this paper, we will show
that Eq. \eqref{GL} is uniformly globally well posed on the
parameter $\nu \ge 0$ in modulation space $M_{2,
1}^s(\mathbb{R}^n),\, n\ge 3,\, s>3$  with the sufficiently small
Cauchy data in $L^2$. As $\nu \to 0$, we prove that the
 solutions of Eq. \eqref{GL} will converge to that of the derivative Schr\"odinger equation. When $n=2$, we
also show  local well-posedness results and inviscid limit in
modulation space $M_{1,1}^s, s>5/2$. The techniques used in this
paper are the anisotropic global smooth effect estimates and maximal
inequality estimates which are independent of parameter $\nu\ge 0$,
those estimates in the case $\nu=0$ were obtained in our earlier
work \cite{WH}, where global well-posedness for equation \eqref{nls}
is showed in $M_{2, 1}^s(\mathbb{R}^n), s\geq 5/2$, for small Cauchy
data.

Finally,  we consider the  quadratic derivative  Ginzburg--Landau
equation:
\begin{align}
u_t- (\nu + {\rm i})\triangle u -
\overrightarrow{\lambda}\cdot\nabla(u^2)=0, \quad
u(0,x)=u_0(x).\label{GL3}
\end{align}
Its limit equation is
\begin{align}
u_t-  {\rm i} \triangle u -
\overrightarrow{\lambda}\cdot\nabla(u^2)=0, \quad
u(0,x)=u_0(x).\label{nls3}
\end{align}
When $ n=1$, Christ in \cite{Ch} showed that for Eq. \eqref{nls3},
the flow map $u_0 \to u$ is not continuous in any Sobolev space
$H^s(\mathbb{R})$ ($s\in \mathbb{R}$) for any short time lifespan
($\|u_0\|_{H^s} \ll 1$ but $\|u(t)\|_{H^s} \gg 1$ for some $t\ll
1$). In \cite{ATA}, Stefanov showed the existence for the weak solutions  in
$H^1$ space with small total disturbance $u_0 \in H^1(\mathbb{R}^1)\cap
L^1(\mathbb{R}^1)\cap \{f: \sup_{x}|\int_{-\infty}^x f(y) dy|\leq
\epsilon\}$.

In this paper, we will show that Eqs. \eqref{GL3} and \eqref{nls3}
are locally well posed in modulation space
$M_{1,1}^{3}(\mathbb{R}^n)$ and the inviscid limit between Eqs.
\eqref{GL3} and \eqref{nls3} also holds in the space
$M_{1,1}^{3}(\mathbb{R}^n)$ for the solutions. From this point of
view, $M_{1,1}^s$ seems to be a proper space to deal with the
solutions of quadratic derivative nonlinear Schr\"odinger equation.

\subsection{Main results}

\begin{thm}\label{thg1}
Let $n\geq 3$. Assume initial data $u_0 \in M_{2,1}^{s}, s>  3$ and
$\|u_0\|_{L^2}\leq \delta$ for some small
$\delta>0$\footnote{$u_0\in M_{2,1}^{s}$ implies that $u_0\in L^2$,
$\delta>0$ may depends on $\|u_0\|_{M_{2,1}^{s}}$. }. Then Eq.
\eqref{GL} has a unique global solution $u_\nu \in C(\mathbb{R}^+,
M_{2,1}^{s})\bigcap X_s$  satisfying
\begin{align}
\|u_\nu\|_{X_s} \leq C \|u_0\|_{M^s_{2,1}} .\nonumber
\end{align}
where $C$ is independent of $\nu$, $X_s$ is defined in \eqref{X}.
\end{thm}

\begin{thm}\label{thg2}
Let $n\geq 3$. Assume initial data $u_0 \in M_{2,1}^{4}$
and $\|u_0\|_{L^2}\leq \delta$ for some small $\delta>0$.
$u_\nu$ is the solution of \eqref{GL}, and let $v$ is the solution
of \eqref{nls} with the same initial data, then for any $T>0$ we
have
\begin{align}
\|u_\nu - v\|_{C(0,T; L^{2})} \lesssim \|u_\nu - v\|_{C(0,T; M^{ }_{2,1})} \lesssim \nu T, \quad \nu \ll 1
.\nonumber
\end{align}
\end{thm}

\begin{thm}
Let $n=1, 2 $. Assume initial data $u_0 \in M_{1,1}^{s}, s> 5/2$ and
$\|u_0\|_{L^{1}}\leq \delta$ for some small
$\delta>0$\footnote{$u_0\in M_{1,1}^{s}$ implies that $u_0\in L^1$,
$\delta>0$ may depends on $\|u_0\|_{M_{1,1}^{s}}$. }. Then Eq.
\eqref{GL} has a unique local solution
$$
u_\nu \in C([0, 1],
M_{2,1}^{s })\bigcap  C([0, 1],
M_{1,1}^{s-1/2})\bigcap X^1_s
$$
 satisfying
$
\|u_\nu\|_{X^1_s} \leq C \|u_0\|_{M_{1,1}^{s}},
$
where $C$ is independent of $\nu$, $X^T$ is defined in \eqref{xt}.
Moreover, if $u_0 \in M_{1,1}^{3}$, then we have
\begin{align}
\|u_\nu - v\|_{C(0,T; L^1)}
\lesssim  \|u_\nu - v\|_{C(0,T; M^{ }_{1,1})}  \lesssim \nu T, \quad
\nu \ll 1 .\nonumber
\end{align}
where $v$ is the solution of the DNLS \eqref{nls} with
the same initial data.
\end{thm}

\begin{thm}
Let $n\in \mathbb{N} $. Assume initial data $u_0 \in M_{1,1}^{s}, s>
3$ and $\|u_0\|_{L^1}\leq \delta$ for some small $\delta>0$. Then
Eq. \eqref{GL3} has a unique solution
$$
 u_\nu \in  C([0, 1],
M_{2,1}^{s })\cap C([0, 1],
M_{1,1}^{s-1/2})\cap {\tilde{X}}^1_s
$$
satisfying $ \|u_\nu\|_{{\tilde{X}}^1_s} \leq C\|u_0\|_{M_{1,1}^{s}}
$, where $C$ is independent of $\nu$, ${\tilde{X}}^T$ is defined in
\eqref{xt1}. Meanwhile,  we have
\begin{align}
 \| u_\nu - v \|_{C(0,T; L^1)} \lesssim \|u_\nu-v\|_{C(0,T; M^{ }_{1,1})}\lesssim \nu T, \quad \nu \ll 1 .\nonumber
\end{align}
where $u_\nu$ and $ v$ is the solution of the DCGL \eqref{GL3} and DNLS \eqref{nls3}
with the same initial data.
\end{thm}

Now we give a brief explanation to the proof of our main results.  We rewrite \eqref{GL} into an integral equation:
\begin{align}
&u= G_{\nu}(t) u_0- \mathscr{A}_\nu
[\overrightarrow{\lambda_1}\cdot\nabla(|u|^2u) +
(\overrightarrow{\lambda_2}\cdot\nabla
u)|u|^2+ \alpha |u|^{2\delta}u],\nonumber\\
& G_{\nu}(t)= \mathscr{F}^{-1}e^{{-\rm i} t|\xi|^2-\nu t|\xi|^2}
\mathscr{F}, \quad \mathscr{A}_\nu f(t, x)=\int_0^t
G_{\nu}(t-\tau)f(\tau, x)d \tau.\nonumber
\end{align}
We will use the smooth effect techniques to prove our result.  Comparing with
the Schr\"odinger equation, the semigroup of Ginzburg-Landau
equation $G_\nu(t)$ dosen't have conjugate symmetry property, this
means
$$G_\nu (t)\neq \overline{G_\nu (-t)},$$
 we can not use
standard $TT^\ast$ argument to get the smooth effect estimates, maximal function estimates and their relations
 with the Strichartz estimates for
$G_\nu (t)$ and  $\mathscr{A}_\nu$. It is known that $TT^*$ method is a basic tool for those estimates
in the case $\nu=0$.

The crucial estimates are  the uniform anisotropic global smooth effect estimates for
semigroup $G_\nu(t)$ and integral operator $\mathscr{A}_\nu$:
\begin{align}
&\|D_{x_i}^{1/2} G_{\nu}(t)\Box_k
u_0\|_{L_{x_i}^{\infty}L_{(x_j)j\neq i }^{2}L_t^2(\mathbb{R}^+
\times \mathbb{R}^n)} \leq C \|\Box_k u_0\|_{2},\quad |k_i|\geq 4,\label{i1}\\
&\|\mathscr{A}_\nu \partial_{x_i}f\|_{L_{x_i}^{\infty}L_{x_j(j\neq
i)}^{2}L_t^2(\mathbb{R}^+ \times \mathbb{R}^n)} \le C
\|f\|_{L_{x_i}^{1}L_{x_j(j\neq i)}^{2}L_t^2(\mathbb{R}^+ \times
\mathbb{R}^n)} \label{i2},
\end{align}
where those estimates in the case $\nu=0$ were established in \cite{14, Li-Po93, WH}.
The main difficulty arises in the fact that the constant $C$ in \eqref{i1} and  \eqref{i2} should be
independent of parameter $\nu \ge 0$. We also need
to show the uniform maximal function estimates for $G_\nu(t)$:
\begin{align}
\|\Box_k G_\nu(t)u_0\|_{L_{x_i}^{2}L_{x_j(j\neq
i)}^{\infty}L_t^{\infty}(\mathbb{R}^+\times \mathbb{R}^n)} \leq C
\langle k_i\rangle^{1/2} \|\Box_k u_0\|_{L^2( \mathbb{R}^n)},\quad n
\geq 3.\label{3.9.1}
\end{align}
In order to show \eqref{3.9.1}, we will use the maximal operator estimates in anisotropic Lebesgue spaces as in \cite{stein1}.
After establishing those uniform estimates, we can use the idea in \cite{WH} to carry out the uniform global well posedness of Eq. \eqref{GL}. The limit behavior can be shown by using the techniques as in \cite{FDF}.

This paper is organized as follows. In Section 2 and Section 3 we
prove the uniform anisotropic global smooth effect estimates, maximal
inequality estimates, Strichartz type estimates for semigroup
$G_\nu(t)$ and integral operator $\mathscr{A}$. In Section 4 we show
the proof of  Theorem 1. In Section 5 we show the proof of inviscid
limit results. In Sections 6 and 7, we show the proof of Theorems 1.3
and   1.4.

\subsection{Notation }
In the sequel $C$ will denote a universal positive constant which
can be different at each appearance. $x\lesssim y $ (for $x$, $y
>0$) means that $x\le Cy$, and $x\sim y$ stands for $x\lesssim y $
and $y\lesssim x$. For any $p\in [1,\infty]$, $p'$ denotes the
conjugate number of $p$, i.e., $1/p+1/p'=1$. Now we introduce the
spaces used in our paper. Let  $\mathscr{S}(\mathbb{R}^n)$ be
Schwartz space. We will use the Lebesgue spaces $L^p:=
L^p(\mathbb{R}^n)$ with the norm $\|\cdot\|_p:=
\|\cdot\|_{L^p(\mathbb{R}^n)}$ and the function spaces $L^q_{t\in I}
L^p_x$ and $L^p_x L^q_{t\in I}$ for which the norms are defined by:
\begin{align}
 \|f\|_{L^q_{t\in I} L^p_x}=\left\|\|f\|_{L^p_x}\right\|_{L^q_t (I)}, \quad  \|f\|_{L^p_x L^q_{t\in I}}=\left\|\|f\|_{L^q_t (I)}\right\|_{L^p_x}. \nonumber
\end{align}
$I$ will be omitted if $I= \mathbb{R}$, i.e., we simply write
$L^q_{t} L^p_x:= L^q_{t\in \mathbb{R}} L^p_x$,~ $L^p_x
L^q_{t}:=L^p_x L^q_{t\in \mathbb{R}}$ and $L^p_x L^q_{T}= L^p_x
L^q_{t\in [0, T]}$,~$L^q_{T}L^p_x =  L^q_{t\in [0, T]} L^p_x$. In
high dimension case, we denote by $L_{x_i}^{p_1}L_{(x_j)j\neq
i}^{p_2}L_t^{p_2}(I\times \mathbb{R}^n)$ the anisotropic Lebesgue
space for which the norm is
\begin{align}
\|f\|_{L_{x_i}^{p_1}L_{x_j(j\neq i)}^{p_2}L_t^{p_2}(I\times
\mathbb{R}^n)}=\Big\|\|f\|_{L_{x_1,\ldots,x_{i-1}, x_{i+1}, \ldots,
x_n}^{p_2}L_t^{p_2}(I\times
\mathbb{R}^{n-1})}\Big\|_{L_{x_i}^{p_1}(\mathbb{R})}.
\end{align}
$D_{x_i}^s = (-\partial_{x_i}^2)^{s/2}=
\mathscr{F}_{\xi_i}^{-1}|\xi_i|^s\mathscr{F}_{x_i}$ denotes the
partial Rieze potential in the $x_i$ direction. $\partial_{x_i}^{-1}
= \mathscr{F}_{\xi_i}^{-1}(i\xi_i)^{-1}\mathscr{F}_{x_i}$. The
homogeneous Sobolev space $\dot{H}^{s}(\mathbb{R}^n)$ is defined by
$(-\Delta)^{-s/2}L^2(\mathbb{R}^n)$.

Modulation spaces were first introduced by Feichtinger \cite{Fei2}. We will use an equivalent norm on the modulation space $M_{2,1}^s$:
\begin{align}
&\|f\|_{M_{2,1}^s}= \sum_{k\in\mathbb{Z}^n} \langle k\rangle^s
\|\mathscr{F}f\|_{L^2(Q_k)},\noindent
\end{align}
where $\langle k\rangle=1+ |k|$, $Q_k= \{\xi: -1/2 \leq \xi_i -k_i
\leq 1/2, i=1, \ldots, n\}$. Let $\{\sigma_k\}_{k\in\mathbb{Z}^n}$ satisfies:
\begin{align}
\left\{
  \begin{array}{ll}
    \sigma_k(\xi)\geq c, & \hbox{$\forall \xi\in Q_k$;} \\
    {\rm supp}\sigma_k\subset \{\xi: |\xi-k|\leq \sqrt{n}\},  \\
    \sum_{k\in\mathbb{Z}^n}\sigma_k(\xi)\equiv 1, & \hbox{$\forall \xi\in \mathbb{R}^n$;} \\
    |D^\alpha \sigma_k(\xi)|\leq C_m, & \hbox{$\forall \xi\in \mathbb{R}^n, |\alpha|\leq m\in \mathbb{N}.$}
  \end{array}\label{denote}
\right.
\end{align}
Denote
\begin{align}
\Upsilon= \{\{\sigma_k\}_{k\in\mathbb{Z}^n}:
\{\sigma_k\}_{k\in\mathbb{Z}^n} \quad
\text{satisfies}\eqref{denote}\,\}.
\end{align}
Let $\{\sigma_k\}_{k\in\mathbb{Z}^n}\in \Upsilon$ be a function
sequence. Then we can define the frequency-uniform decomposition
operators $\Box_k$ as:
\begin{align}
\Box_k := \mathscr{F}^{-1}\sigma_k \mathscr{F}, \quad
k\in\mathbb{Z}^n,
\end{align}
and we have
\begin{align}
&\|f\|_{M_{2,1}^s}\sim\sum_{k\in\mathbb{Z}^n} \langle k\rangle^s
\|\Box_k f\|_{L^2(\mathbb{R}^n)}.
\end{align}
Using the operators $\Box_k$, we can equivalently define the modulation space $M_{1,1}^s$ in the following way:
\begin{align}
&\|f\|_{M_{1,1}^s}= \sum_{k\in\mathbb{Z}^n} \langle k\rangle^s
\|\Box_k f\|_{L^1(\mathbb{R}^n)}.
\end{align}

For simplicity, we use function space
$l_{\Box}^{1,s}(L^{\nu}(\mathbb{R}^+; L^r(\mathbb{R}^n)))$ which
contains all of the functions $f(t, x)$ so that the following norm
is finite:
\begin{align}
\sum_{k\in\mathbb{Z}^n} \langle k\rangle^s \|\Box_k f\|_{L_t^\nu
L_x^r(\mathbb{R}^+\times\mathbb{R}^n)}.
\end{align}

\section{Anisotropic Global smooth effect with $\Box_k$-decomposition}

In this section, we always denote
\begin{align}
G_{\nu}(t)= \mathscr{F}^{-1}e^{{-\rm i} t|\xi|^2}\cdot
e^{-\nu t |\xi|^2}\mathscr{F},\quad \mathscr{A}f(t,
x)=\int_0^t G_{\nu}(t-\tau)f(\tau, x)d \tau.
\end{align}

For convenience, we will use the following function sequence
$\{\sigma_k\}_{k\in \mathbb{Z}^n}$:
\begin{lem}\label{ll}
Let $\eta_k : \mathbb{R}\rightarrow [0, 1] (k\in \mathbb{Z})$ be a
smooth-function sequence satisfying \eqref{denote}. Denote
\begin{align}
\sigma_k(\xi) := \eta_{k_1}(\xi_1)\ldots  \eta_{k_n}(\xi_n), \quad
k=(k_1, \ldots, k_n) \label{2.2}.
\end{align}
Then we have $\{\sigma_k\}_{k\in \mathbb{Z}^n} \in \Upsilon$.
\end{lem}

\begin{lem}{\rm (\cite{WH})}\label{lem2.4}
We have for any $\delta\in \mathbb{R}$ and $k=(k_1, \ldots, k_n)\in
\mathbb{Z}^n$ with $|k_i| \geq 4$,
\begin{align}
\|\Box_k D_{x_i}^\sigma u\|_{L_{x_1}^{p_1}L_{x_2, \ldots
x_n}^{p_2}L_t^{p_2}(\mathbb{R}_+\times \mathbb{R}^n)} \lesssim
\langle k_i \rangle^{\sigma}\|\Box_k u\|_{L_{x_1}^{p_1}L_{x_2,
\ldots x_n}^{p_2}L_t^{p_2}(\mathbb{R}_+\times \mathbb{R}^n)}.
\end{align}
Replacing $D_{x_i}^\delta$ by $\partial_{x_i}^\delta (\delta\in
\mathbb{N})$, the above inequality holds for all $k\in
\mathbb{Z}^n$.
\end{lem}

In order to obtain global smooth-effect estimates, we
need the following Lemma in the case $n=1$:

\begin{lem}\label{p1.3}
Let $n=1, |k|\geq 4$. Then  there exists $C>0$, which is
independent of $\nu>0$ such that
\begin{align}
\|\mathscr{F}_{\xi}^{-1}e^{-i t\xi^2}e^{-\nu
|t|\xi^2}\mathscr{F}_{x}\Box_k
\phi\|_{L_{x}^{\infty}L_t^2(\mathbb{R}\times \mathbb{R} )}\leq
C\langle k\rangle^{-1/2}\| \Box_k \phi\|_{L_2(\mathbb{R})}.
\end{align}
\end{lem}
\noindent{\bf Proof.} We may assume that $k\ge 4$. By  changing the variable, we have:
\begin{align}
 \mathscr{F}_{\xi}^{-1} & e^{-i t\xi^2}e^{-\nu
|t|\xi^2} \mathscr{F}_{x}\Box_k\phi\nonumber\\
&=\int_{\mathbb{R}} e^{{-\rm i} t\xi^2}e^{-\nu
|t|\xi^2}\eta_k(\xi)\hat{\phi}(\xi) e^{{\rm i}x \xi}d \xi \nonumber\\
&=\int_{\mathbb{R}} e^{{-\rm i} t\eta}e^{-\nu
|t|\eta}\eta_k(\sqrt{\eta})\hat{\phi}(\sqrt{\eta}) e^{{\rm i}x
\sqrt{\eta}}\frac{1}{2\sqrt{\eta}}d \eta,
\end{align}
where $\eta_k$  was defined in Lemma \ref{ll}.

 From Plancherel's
identity, Fubini theorem and Young's inequality we have
\begin{align}
&\|\mathscr{F}_{\xi}^{-1}e^{-i t\xi^2}e^{-\nu
|t|\xi^2} \mathscr{F}_{x}\Box_k
\phi\|_{L_{x}^{\infty}L_t^2(\mathbb{R}\times \mathbb{R} )}\nonumber\\
&\leq \Big\|\int_\mathbb{R} \mathscr{F}_t(e^{-\nu
|t|\eta})(\tau+\eta)\Big|\eta_k(\sqrt{\eta})\hat{\phi}(\sqrt{\eta})
e^{{\rm i}x
\sqrt{\eta}}\Big|\frac{1}{2\sqrt{\eta}}d \eta \Big\|_{L_x^\infty L_{\tau}^2}\nonumber\\
&=\Big\|\int_\mathbb{R} \Big(\frac{1}{\nu
\eta}\frac{1}{1+\Big|\frac{\tau+\eta}{\nu
\eta}\Big|^2}\Big)\Big|\eta_k(\sqrt{\eta})\hat{\phi}(\sqrt{\eta})
\Big|\frac{1}{2\sqrt{\eta}}d \eta \Big\|_{L_{\tau}^2}\nonumber\\
 &\leq \Big\|\int_\mathbb{R} \Big(\frac{1}{\nu (\langle
k\rangle-1/2)^2}\frac{1}{1+\Big|\frac{\tau+\eta}{\nu (\langle
k\rangle+1/2)^2}\Big|^2}\Big)\Big|\eta_k(\sqrt{\eta})\hat{\phi}(\sqrt{\eta})
\Big|\frac{1}{2\sqrt{\eta}}d \eta \Big\|_{L_{\tau}^2}\label{ti}\\
 &\lesssim \Big\|\eta_k(\sqrt{\eta})\hat{\phi}(\sqrt{\tau})
\frac{1}{2\sqrt{\tau}}\Big\|_{L_{\tau}^2} \Big\|
\mathscr{F}_t(e^{-\nu |t|(\langle
k\rangle+1/2)^2})(\tau)\Big\|_{L_{\tau}^1}\nonumber\\
&\leq C \langle k\rangle^{-1/2}\|\Box_k \phi\|_{L_2},\nonumber
\end{align}
where we have used
\begin{align}
\Big\| \mathscr{F}_t(e^{-\nu
|t|\eta})(\tau)\Big\|_{L_{\tau}^1}=\int
\frac{1}{1+\tau^2}d \tau \leq C,
\end{align}
where $C$ is independent of $\nu$ and $\eta$. In equation \eqref{ti},
$|k|\geq 4$ is necessary. $\hfill\Box$

\begin{lem} \label{prop4}
\begin{align}
\sup_{s\geq 0}\Big\|\mathscr{F}^{-1}_{\tau, \xi}
\frac{\xi}{\xi^2-i\nu (\xi^2 + s)+\tau } \mathscr{F}_{t,
x}f\Big\|_{L_{x}^{\infty}L_t^2(\mathbb{R}^{1+1})} \leq C
\|f\|_{L_{x}^{1}L_t^2(\mathbb{R}^{1+1})}.\label{1.2}
\end{align}
The constant $C$ in \eqref{1.2} is independent of $\nu>0$
.\label{lem1}
\end{lem}
\noindent {\bf Proof.}  For convenience, we denote by
$\mathscr{F}_{t,x}$, $\mathscr{F}_t$,
 $\mathscr{F} $ the Fourier transforms on $(t,x)$, $t$, $x$,
respectively. From Plancherel's identity, we have :
\begin{align}
&\quad \|\mathscr{F}^{-1}_{\xi}\frac{\xi}{\tau+\xi^2- i\nu (\xi^2+
s)}
\mathscr{F}_{t,x} f\|_{L_{\tau}^2}\nonumber\\
&=\Big\|\int_{\mathbb{R}}\int_{\mathbb{R}} e^{i(x-y)\xi}
\frac{\xi}{\tau+\xi^2 - i\nu (\xi^2+s)} (\mathscr{F}_t f)(\tau, y)
d\xi dy
\Big\|_{L_{\tau}^2}\nonumber\\
&:=\Big\|\int_\mathbb{R} K(\tau, x-y) (\mathscr{F}_t f) (\tau, y)dy
\Big\|_{L_{\tau}^2} \label{}.
\end{align}
where the integral
\begin{align}
K(\tau, z)&=\int_\mathbb{R} e^{iz\xi} \frac{\xi}{\tau+\xi^2 - i\nu
(\xi^2 + s)}d\xi \label{}
\end{align}
is taken in the P.V. meaning. Now we only need to show that
\begin{align}
\sup_{s\geq 0}\|K(\tau, z)\|_{L^\infty_{\tau, z}} \lesssim 1,
\label{1.5}
\end{align}
We only consider case $\tau<0$, (the case $\tau\geq 0$ do not
contain singular point, so it is easy to handle). For $\tau<0$, we
have
\begin{align}
K(\tau, z)&=\int_\mathbb{R} e^{iz\sqrt{-\tau}\eta}
\frac{\eta}{1-\eta^2+ i\nu (\eta^2+s_1)} d\eta \nonumber\\
&=\int_\mathbb{R} e^{iz\sqrt{-\tau}\eta}\  \frac{\eta[1-\eta^2- i\nu
(\eta^2+s_1)]}{\nu^2 (\eta^2+s_1)^2+(1-\eta^2)^2} d\eta,
\end{align}
where $s_1=-s/\tau>0$.

Since when $s_1\rightarrow 0, \nu\rightarrow 0, \eta\rightarrow 1$,
$K(\tau, z)$ is difficult to handle, we will divide $\eta$ into
different cases: Let $\psi_1, \psi_2, \psi_3 \in
C_0^{\infty}(\mathbb{R})$ satisfy ${\rm supp} \psi_1\subseteq
\{\eta: |\eta|\geq 3/2\}, \psi_1(-\cdot)=\psi_1(\cdot)$, ${\rm
supp}\psi_2 \subseteq (-2, 1/2]$, ${\rm supp}\psi_3 \subseteq (0,
2)$, $\sum_{i=1}^3\psi_i=1$. Define

\begin{align}
K_i(\tau, z) &=\int_\mathbb{R} e^{iz\sqrt{-\tau}\eta}\
\frac{\eta[1-\eta^2-i\nu (\eta^2+s_1)]\psi_i(\eta)}{\nu^2
(\eta^2+s_1)^2+(1-\eta^2)^2} d\eta,\ \ \ \ i=1,2,3.
\end{align}
\begin{align}
K_1(\tau, z) =&\int_{|\eta|\geq 3/2} e^{iz\sqrt{-\tau}\eta}\
\frac{\eta}{\nu^2 (\eta^2+s_1)^2+(1-\eta^2)^2}d\eta
\nonumber\\
&-(1+i \nu)\int_{|\eta|\geq 3/2} e^{iz\sqrt{-\tau}\eta}\
\frac{\eta^3}{\nu^2 (\eta^2+s_1)^2+(1-\eta^2)^2}d\eta
\nonumber\\
&-i \int_{|\eta|\geq 3/2} e^{iz\sqrt{-\tau}\eta}\ \frac{\nu\eta
s_1}{\nu^2
(\eta^2+s_1)^2+(1-\eta^2)^2}d\eta\nonumber\\
 &:= K_1^1(\tau, z)+K_1^2(\tau,
z)+K_1^3(\tau, z).
\end{align}
 It is easy to see
\begin{align}
|K_1^1(\tau, z)| \lesssim \int_{|\eta|\geq 3/2} \frac{1}{\eta^3} d
\eta\leq C.
\end{align}
From variable changing, we have:
\begin{align}
&|K_1^3(\tau, z)|\leq 2 \int_{\eta\geq 3/2} \frac{\nu s_1 \eta
}{(\nu s_1)^2+ \eta^4}d \eta \leq \arctan \eta \Big|_{\eta/\nu
s_1}^{\infty} \leq C,\label{1.10}
\end{align}
\begin{align}
|K_1^2(\tau, z)| \lesssim \Big|\int_{|\eta|\geq 3/2}
e^{iz\sqrt{-\tau}\eta}\ \frac{\eta^3}{\nu^2
(\eta^2+s_1)^2+(1-\eta^2)^2} d\eta\Big|.
\end{align}
We derive $K_1^2(\tau, z)$ into two parts $I$, $II$, from
\eqref{1.10} we have:
\begin{align}
& I=\Big|\int_{ 3/2 \leq |\eta| \leq 10\sqrt{\nu s_1}}
e^{iz\sqrt{-\tau}\eta}\ \frac{\eta^3}{\nu^2
(\eta^2+s_1)^2+(1-\eta^2)^2} d\eta\Big|\nonumber\\
&  \lesssim \Big|\int_{ 3/2 \leq |\eta| \leq 10\sqrt{\nu
s_1}}\frac{\nu s_1 \eta }{(\nu s_1)^2+ \eta^4} \leq \arctan \eta
\Big|_{\eta/\nu s_1}^{\infty}\leq C.
\end{align}
From variable changing, we have:
\begin{align}
& |II|=  \Big|\int_{   \eta \geq {10\sqrt{\nu s_1}}}
\sin(z\sqrt{-\tau}\eta)
\frac{1}{(1+\nu^2)\eta + \frac{\nu^2 s_1^2 +1}{\eta^3} + \frac{2\nu^2 s_1-2}{\eta}} d\eta\Big|\nonumber\\
&=\Big|\int_{   \eta \geq {10z\sqrt{-\tau\nu s_1}}} \sin(\eta)
\frac{1}{(1+\nu^2)\eta + \frac{(z\sqrt{-\tau})^4(\nu^2 s_1^2 +1)}
{\eta^3} + \frac{(2\nu^2 s_1-2)(z\sqrt{-\tau})^2 }{\eta}}
d\eta\Big|\label{1.13}.
\end{align}
Now we prove \eqref{1.13}is bounded. Write
$$F(\eta):=\frac{1}{(1+\nu^2)\eta +
\frac{(z\sqrt{-\tau})^4(\nu^2 s_1^2 +1)}{\eta^3} + \frac{(2\nu^2
s_1-2)(z\sqrt{-\tau})^2 }{\eta} }.$$ For any $\epsilon >0$,  when
$A'> A
> 1/(1+\nu^2)\epsilon$, we have:
$$F(A)\leq \epsilon, \quad F(A)'\leq \epsilon.$$
Notice that $F(\eta)$ is monotonous decreasing when $ \eta \geq
{10\sqrt{\nu s_1}}$ and for any $\eta\in [A, A']$,
$\int_A^{A'}sin(\eta)d \eta\leq C$. So from the second
integral-mean-value theorem, we have:
\begin{align}
\Big|\int_A^{A'} \sin(\eta) \frac{1}{(1+\nu^2)\eta +
\frac{(z\sqrt{-\tau})^4(\nu^2 s_1^2 +1)} {\eta^3} + \frac{2\nu^2
s_1(z\sqrt{-\tau})^2 }{\eta}} d\eta\Big|\leq C\epsilon.
\end{align}
Then from the Cauchy convergence theorem, we can get \eqref{1.13} is
bounded. So
$$
\|K_1(\tau, z)\|_{L^\infty_{\tau, z}} \lesssim 1.
$$
Notice that $\nu^2 (\eta^2+s_1)^2+(1-\eta|\eta|)^2\geq 3/4$, when
$\eta \in (-3/2, 1/2]$, so it is easy to estimate $K_2(\tau, z)$:
\begin{align}
\|K_2(\tau, z)\|_{L^\infty_{\tau, z}} &\leq \int_\mathbb{R}
\Big|e^{iz\sqrt{-\tau}\eta}\ \frac{\eta[1-\eta^2- i\nu
(\eta^2+s_1)]\psi_2(\eta)}{\nu^2 (\eta^2+s_1)^2+(1-\eta^2)^2}\Big|
d\eta  \nonumber\\
& \lesssim \Big| \frac{\eta[1-\eta^2- i\nu
(\eta^2+s_1)]\psi_2(\eta)}{\nu^2 (\eta^2+s_1)^2+(1-\eta^2)^2}\Big|
\leq C.
\end{align}
\begin{align}
K_3(\tau, z)&=\int_\mathbb{R} e^{iz\sqrt{-\tau}\eta}
\frac{\eta\psi_3(\eta)}{1-\eta^2+ i\nu (\eta^2+ s_1)} d\eta \nonumber\\
&=\int_\mathbb{R}
e^{iz\sqrt{-\tau}\eta}\frac{\eta\psi_3(\eta)}{1-\eta^2} d \eta -
\int_\mathbb{R} e^{iz\sqrt{-\tau}\eta} \frac{i \nu \eta
(\eta^2+s_1)\psi_3(\eta)}{[1-\eta^2+ i\nu (\eta^2 + s_1)](1-\eta^2)}
d\eta
\nonumber\\
&:= K_3^1(\tau, z) + K_3^2(\tau, z),\label{2.29.0}
\end{align}
\begin{align}
|K_3^1(\tau, z)| &= \Big|\int_\mathbb{R} e^{iz\sqrt{-\tau}\eta}
\frac{1}{1-\eta}\cdot\frac{\eta\psi_3(\eta)}{1+\eta} \ d \eta\Big|\
&\sim\Big|\Big(sgn \ast \mathscr{F}^{-1}\Big[
\frac{\eta\psi_3(\eta)}{1+\eta} \Big]\Big)(\sqrt{-\tau
z})\Big|,\label{2.29.1}
\end{align}
\begin{align}
|K_3^2(\tau, z)| &=\Big|\int_\mathbb{R} e^{iz\sqrt{-\tau}\eta}
\frac{1}{1-\eta}\cdot\frac{i \nu \eta(\eta^2+s_1)[1-\eta^2-i\nu
(\eta^2+s_1)]\psi_3(\eta)}{[(1-\eta^2)^2+ \nu^2 (\eta^2+
s_1)^2](1+\eta)} d\eta\Big|
\nonumber\\
&\sim\Big|\Big(sgn \ast \mathscr{F}^{-1}\Big[ \frac{i \nu
\eta(\eta^2+s_1)[1-\eta^2-i\nu
(\eta^2+s_1)]\psi_3(\eta)}{[(1-\eta^2)^2+ \nu^2 (\eta^2
+s_1)^2](1+\eta)} \Big]\Big)(\sqrt{-\tau z})\Big|.\label{2.29.2}
\end{align}
Noticing that
$$
\Big|\frac{i \nu \eta(\eta^2+s_1)[1-\eta^2-i\nu
(\eta^2+ s_1)]\psi_3(\eta)}{[(1-\eta^2)^2+ \nu^2
(\eta^2+s_1)^2](1+\eta)}\Big| \leq C,
$$
from \eqref{2.29.0}--\eqref{2.29.2}, we obtain that
\begin{align}
\|K_3(\tau, z)\|_{L^\infty_{\tau, z}} \leq  &
\Big\|\mathscr{F}^{-1}\Big[ \frac{\eta\psi_3(\eta)}{1+\eta}
\Big]\Big\|_{L_1} \nonumber\\
& +\Big\|\mathscr{F}^{-1}\Big[ \frac{i \nu
\eta(\eta^2+s_1)[1-\eta^2-i\nu
(\eta^2+s_1)]\psi_3(\eta)}{[(1-\eta^2)^2+ \nu^2 (\eta^2
+s_1)^2](1+\eta)} \Big]\Big\|_{L_1}\lesssim 1.
\end{align}
It follows that \eqref{1.5} holds.

\begin{prop}\label{p1.4}
For any $i=1, \ldots, n$, $|k_i|\geq 4$,  we have
\begin{align}
 &\|D_{x_i}^{1/2} \Box_k
G_{\nu}(t)u_0\|_{L_{x_i}^{\infty}L_{(x_j)j\neq i
}^{2}L_t^2(\mathbb{R}^+\times \mathbb{R}^n)} \lesssim \|\Box_k
u_0\|_{2}.\label{1.7}
\end{align}
\end{prop}
\noindent{\bf{Proof.}}
 As in the proof
of Lemma \ref{p1.3}, we only need to prove that
\begin{align}
&\|D_{x_i}^{1/2} \Box_k
G_{\nu}'(t)u_0\|_{L_{x_i}^{\infty}L_{(x_j)j\neq i
}^{2}L_t^2(\mathbb{R} \times \mathbb{R}^n)} \lesssim \|\Box_k
u_0\|_{2}, \label{1.2.1}
\end{align}
where $G_{\nu}(t)'= \mathscr{F}^{-1}e^{{-\rm i}
t|\xi|^2-\nu |t| |\xi|^2} \mathscr{F}$.
It suffices to show the case $i=1$.  By Plancherel's identity
and Minkowski's inequality,
\begin{align}
&\| D_{x_1}^{1/2}G_\nu(t)'\Box_k u_0\|_{L_{x_1}^{\infty}L_{x_2
\ldots
x_n}^{2}L_t^2(\mathbb{R}^{1+n})} \nonumber\\
&\leq \|D_{x_1}^{1/2}\mathscr{F}_{\xi_1}^{-1}e^{-i t\xi_1^2}e^{-
\nu|t|\xi_1^2}\mathscr{F}_{x_1}(\mathscr{F}_{x_2, \ldots, x_n}\Box_k
u_0)\|_{L_{x_1}^{\infty}L_{\xi_2 \ldots
\xi_n}^{2}L_t^2(\mathbb{R}^{1+n})}\nonumber\\
&\leq \|D_{x_1}^{1/2}\mathscr{F}_{\xi_1}^{-1}e^{-i t\xi_1^2}e^{-
\nu|t|\xi_1^2}\mathscr{F}_{x_1}(\mathscr{F}_{x_2, \ldots, x_n}\Box_k
u_0)\|_{L_{\xi_2 \ldots
\xi_n}^{2}L_{x_1}^{\infty}L_t^2(\mathbb{R}^{1+n})}.\nonumber
\end{align}
 In view of Lemma \ref{p1.3} in one spatial dimension, using Plancherel's
identity, we immediately obtain \eqref{1.2.1}. $\hfill\Box$

\begin{prop}\label{4.12}
For any  $k=(k_1, \ldots, k_n)\in\mathbb{Z}^n$, $|k_i| \ge 4$,  we have
\begin{align}
\|\partial_{x_i} \Box_k \mathscr{A}_\nu
f\|_{L_t^{\infty}L_x^2(\mathbb{R}_+\times \mathbb{R}^n)} \lesssim
\langle k_i\rangle^{1/2}\|\Box_k f\|_{L_{x_i}^{1}L_{(x_j)j\neq i
}^{2}L_t^2(\mathbb{R}_+\times \mathbb{R}^n)}.\label{p1.5}
\end{align}
\end{prop}
\noindent{\bf{Proof.}} For $|k_i|\ge 4$, from Proposition \ref{p1.4},
\eqref{1.7} has the following dual estimate:
\begin{align}
\Big\|\Box_k D_{x_i}^{1/2}\int^t_{0}G_{\nu}(t-\tau)f(\tau)d
\tau\Big\|_{ L_t^\infty L_x^2( \mathbb{R}^n)}\lesssim \|\Box_k
f\|_{L_{x_i}^{1}L_{(x_j)j\neq i }^{2}L_t^2(\mathbb{R}_+\times
\mathbb{R}^n)}.\label{2.20}
\end{align}
Then from Lemma \ref{lem2.4}, which implies \eqref{p1.5} holds, as
desired.
 $\hfill\Box$

\begin{prop} \label{prop2.1}
For any $i=1,\ldots , n$ and $k=(k_1, \ldots, k_n)$,  there exist
$C>0$, which are independent of $\nu>0$ such that
\begin{align}
&\|\Box_k \mathscr{A}_\nu
\partial_{x_i}f\|_{L_{x_i}^{\infty}L_{x_j(j\neq
i)}^{2}L_t^2(\mathbb{R}^+\times \mathbb{R}^n)} \le C \|\Box_k
f\|_{L_{x_i}^{1}L_{x_j(j\neq i)}^{2}L_t^2(\mathbb{R}^+\times
\mathbb{R}^n)}.\label{1.22}
\end{align}
\end{prop}
\noindent {\bf Proof.} In order to prove \eqref{1.22}, assume $f(t,
x)=0$ for $t<0$, so we only need to prove
\begin{align}
\|\Box_k
\mathscr{A}_\nu\partial_{x_1}f\|_{L_{x_1}^{\infty}L_{x_j(j\neq
1)}^{2}L_t^2(\mathbb{R} \times \mathbb{R}^n)} \le C \|\Box_k
f\|_{L_{x_1}^{1}L_{x_j(j\neq 1)}^{2}L_t^2(\mathbb{R}\times
\mathbb{R}^n)} \label{}.
\end{align}
  We have
\begin{align}
\partial_{x_1}\mathscr{A}_\nu f = C \mathscr{F}^{-1}_{\tau, \xi}
\frac{\xi_1}{\xi_1^2-i\nu (\xi_1^2 +|\bar{\xi}|^2)+\tau +
|\bar{\xi}|^2} \mathscr{F}_{t, x} f, \label{a4}
\end{align}
where we assume that the right hand side of \eqref{a4} is zero as
$t=0$. It follows that
\begin{align}
&\|\partial_{x_1}\mathscr{A}_\nu f\|_{L_{x_1}^{\infty}L_{x_2 \ldots x_n}^{2}L_t^2(\mathbb{R}^{1+n})}\nonumber\\
&\leq \Big\|\mathscr{F}^{-1}_{\xi_1} \frac{\xi_1}{\xi_1^2-i\nu
(\xi_1^2 +|\bar{\xi}|^2)+(\tau + |\bar{\xi}|^2)} \mathscr{F}_{t,
x}f\Big\|_{L_{\xi_2\ldots
\xi_n}^{2}L_{x_1}^{\infty}L_\tau^2(\mathbb{R}^{1+n})}.\label{a1}
\end{align}
Now changing the variable $\tau+ |\bar{\xi}|^2\rightarrow \mu $, we
have
\begin{align}
&\Big\|\mathscr{F}^{-1}_{\xi_1} \frac{\xi_1}{\xi_1^2-i\nu (\xi_1^2
+|\bar{\xi}|^2)+(\tau + |\bar{\xi}|^2)} \mathscr{F}_{t,
x}f\Big\|_{L_{\xi_2\ldots
\xi_n}^{2}L_{x_1}^{\infty}L_\tau^2(\mathbb{R}^{1+n})}\nonumber\\
&\leq \sup_{s\geq 0}\Big\|\mathscr{F}^{-1}_{\xi_1}
\frac{\xi_1}{\xi_1^2-i\nu (\xi_1^2 + s)+\mu} \mathscr{F}_{t,
x_1}(e^{{\rm i} t|\bar{\xi}|^2}\mathscr{F}_{x_2, \ldots,
x_n}f)\Big\|_{L_{\xi_2\ldots
\xi_n}^{2}L_{x_1}^{\infty}L_\mu^2(\mathbb{R}^{1+n})}.\label{a2}
\end{align}
From the uniform smooth effect estimate
as in Lemma \ref{lem1},
\begin{align}
\sup_{s\geq 0}\Big\|\mathscr{F}^{-1}_{\tau, \xi}
\frac{\xi}{\xi^2-i\nu (\xi^2 + s)+\tau } \mathscr{F}_{t,
x}f\Big\|_{L_{x}^{\infty}L_t^2(\mathbb{R}^{1+1})} \leq C
\|f\|_{L_{x}^{1}L_t^2(\mathbb{R}^{1+1})}\label{a3}.
\end{align}
From \eqref{a1}, \eqref{a2} and \eqref{a3}, we have that
\begin{align}
\|\partial_{x_1}\mathscr{A}_\nu f\|_{L_{x_1}^{\infty}L_{x_2 \ldots
x_n}^{2}L_t^2(\mathbb{R}^{1+n})} \lesssim \big\|e^{{\rm i}
t|\bar{\xi}|^2}\mathscr{F}_{x_2, \ldots, x_n}f \big\|_{L_{\xi_2,
\ldots ,\xi_n}^{2}L_{x_1}^{1}L_t^2(\mathbb{R}^{1+n})}.
\end{align}
Using Minkowski's inequality and Plancherel's identity, we
immediately have
\begin{align}
\|\partial_{x_1}\mathscr{A}_\nu f\|_{L_{x_1}^{\infty}L_{x_2 \ldots
x_n}^{2}L_t^2(\mathbb{R}^{1+n})} \lesssim \|f\|_{L_{x_1}^{1}L_{x_2,
\ldots ,x_n}^{2}L_t^2(\mathbb{R}^{1+n})}.
\end{align}
Other cases can be shown in a similar way.

Generally, the right hand side in \eqref{a4} is not equal to zero for $t=0$:
\begin{align}
&{\rm i}\mathscr{F}^{-1}_{\tau, \xi} \frac{1}{\tau + |\xi|^2-i\nu
|\xi|^2}
\mathscr{F}_{t, x}f\Big|_{t=0}\nonumber\\
&=-{\rm i}\mathscr{F}^{-1}_{\tau, \xi} \frac{\tau + |\xi|^2}{(\tau +
|\xi|^2)^2+(\nu |\xi|^2)^2} \mathscr{F}_{t, x}f\Big|_{t=0} +
\mathscr{F}^{-1}_{\tau, \xi} \frac{\nu|\xi|^2}{(\tau +
|\xi|^2)^2+(\nu |\xi|^2)^2} \mathscr{F}_{t,
x}f\Big|_{t=0}\nonumber\\
&:=u_1(0, x)+u_2(0, x).
\end{align}
Noticing that $\mathscr{F}_t(e^{-\nu |t|\eta})(\tau)=\frac{1}{\nu
\eta}\frac{1}{1+\big|\frac{\tau}{\nu \eta}\big|^2} $ and changing the
variable,   we have
\begin{align}
u_2(0, x)&= C\int_{\mathbb{R}^n}e^{{\rm i}x\xi}\int_\mathbb{R}
\frac{\nu|\xi|^2}{(\tau + |\xi|^2)^2+(\nu
|\xi|^2)^2}\int_\mathbb{R}e^{{-\rm i}\tau s}\hat{f}(s, \xi)ds d\tau
d\xi\nonumber\\
&= C\int_{\mathbb{R}^n}e^{{\rm i}x\xi}\int_\mathbb{R}\hat{f}(s, \xi)
\int_\mathbb{R}\frac{e^{{-\rm i}\tau s}\nu|\xi|^2}{(\tau +
|\xi|^2)^2+(\nu |\xi|^2)^2} d\tau ds
d\xi\nonumber\\
&= C\int_{\mathbb{R}^n}e^{{\rm i}x\xi}\int_\mathbb{R}\hat{f}(s, \xi)
 \ e^{{\rm i}s|\xi|^2}\int_\mathbb{R}\frac{e^{{-\rm i}\tau
s}\nu|\xi|^2}{\tau^2+(\nu |\xi|^2)^2} d\tau ds
d\xi\nonumber\\
&= C\int_{\mathbb{R}^n}e^{{\rm i}x\xi}\int_\mathbb{R}\hat{f}(s, \xi)
 \ e^{{\rm i}s|\xi|^2} e^{-\nu|s||\xi|^2} ds
d\xi\nonumber\\
&=C\int_\mathbb{R} G_\nu(-s)f(s, x) ds.\label{u2}
\end{align}
Then from \eqref{1.7}, \eqref{u2} and \eqref{2.20} we have
\begin{align}
\|\Box_k
\partial_{x_i}G_\nu(t)u_2(0, x)\|_{L_{x_i}^{\infty}L_{x_j(j\neq
i)}^{2}L_t^2(\mathbb{R}^+\times \mathbb{R}^n)} &\leq
\|D_{x_i}^{1/2}\Box_k u_2(0, x)\|_{L^2}\nonumber\\
&\leq \|D_x^{1/2}\Box_k \int_\mathbb{R} G_\nu(-s)f(s, x)
ds\|_{L^2}\nonumber\\
 &\leq \|\Box_k f\|_{L_{x_i}^{1}L_{x_j(j\neq
i)}^{2}L_t^2(\mathbb{R} \times \mathbb{R}^n)}.\label{2.38}
\end{align}
 Noticing that
\begin{align}
\mathscr{F}_t(\partial_te^{-
|t|}(\nu|\xi|^2\cdot))(\tau)=\frac{1}{\nu|\xi|^2}\mathscr{F}_t
(\partial_te^{-
|t|})(\frac{\tau}{\nu|\xi|^2})=\frac{1}{\nu|\xi|^2}\frac{\tau}
{\nu|\xi|^2}\frac{1}{1+(\tau/\nu|\xi|^2)^2}
=\frac{\tau}{\tau^2+(\nu|\xi|^2)^2}, \nonumber
\end{align}
similar to \eqref{u2}, we have
\begin{align}
u_1(0, x)&= C\int_{\mathbb{R}^n}e^{{\rm i}x\xi}\int_\mathbb{R}
\frac{\tau+|\xi|^2}{(\tau + |\xi|^2)^2+(\nu
|\xi|^2)^2}\int_\mathbb{R}e^{{-\rm i}\tau s}\hat{f}(s, \xi)ds d\tau
d\xi\nonumber\\
&= C\int_{\mathbb{R}^n}e^{{\rm i}x\xi}\int_\mathbb{R}\hat{f}(s,
\xi)e^{{\rm i}s|\xi|^2} \int_\mathbb{R}\frac{e^{{-\rm i}\tau
s}\tau}{\tau^2+(\nu |\xi|^2)^2} d\tau ds d\xi\nonumber\\
&= C\int_{\mathbb{R}^n}e^{{\rm i}x\xi}\int_\mathbb{R}\hat{f}(s, \xi)
 \ e^{{\rm i}s|\xi|^2} {\rm sgn}(s)e^{-\nu|s||\xi|^2} ds
d\xi\nonumber\\
&=C\int_\mathbb{R} G_\nu(-s) {\rm sgn}(s)f(s, x) ds.\label{u1}
\end{align}

Then from \eqref{1.7}, \eqref{u1} and \eqref{2.20} we have
\begin{align}
\|\Box_k
\partial_{x_i}G_\nu(t)u_1(0, x)\|_{L_{x_i}^{\infty}L_{x_j(j\neq
i)}^{2}L_t^2(\mathbb{R}^+\times \mathbb{R}^n)}
&\leq\|D_{x_i}^{1/2}\Box_k u_1(0, x)\|_{L^2}\nonumber\\
&\leq \Big \|D_{x_i}^{1/2}\Box_k
\int_\mathbb{R} G_\nu(-s) {\rm sgn}(s)f(s, x) ds\Big\|_{L^2}\nonumber\\
 &\leq \|\Box_k [{\rm sgn}(s)f(s, \xi)]\|_{L_{x_i}^{1}L_{x_j(j\neq
i)}^{2}L_s^2(\mathbb{R} \times \mathbb{R}^n)}\nonumber\\
&\leq \|\Box_k f(s, \xi)\|_{L_{x_i}^{1}L_{x_j(j\neq
i)}^{2}L_s^2(\mathbb{R} \times \mathbb{R}^n)}.\label{2.40}
\end{align}
Collecting \eqref{2.40}, \eqref{2.38}, we can obtain the result, as
desired. $\hfill\Box$

\section{Other estimates with $\Box_k$-decomposition}
In this section, we consider the Strichartz estimates, the maximal
function estimates and derivative interaction estimates
for the solutions of Ginzburg-Laundau equation by using the
frequency-uniform decomposition operators.

\vspace{10pt}

Using Lemma \ref{lem2.2} and the property of frequency-uniform
decomposition operators (cf.\cite{WGZ}), we can establish the
following Strichartz estimates in a class of function spaces by
using the frequency-uniform decomposition operators.

\begin{prop}\label{prop3.3}
Let $2 \leq r  <\infty$, $q > \nu \geq 2\vee \nu(r)$, then we have
\begin{align}
&\|G_\nu(t)f\|_{l_{\Box}^1(L^{\nu}(\mathbb{R}_+; L^r(\mathbb{R}^n)))}\leq C \|f\|_{M_{2, 1}(\mathbb{R}^n)},\\
&\|\mathscr{A}_\nu f\|_{l_{\Box}^1(L^{\nu}(\mathbb{R}_+;
L^r(\mathbb{R}^n)))\bigcap l_{\Box}^1(L^{\infty}(\mathbb{R}_+;
L^2(\mathbb{R}^n)))}\leq C \|f\|_{l_{\Box}^1(L^{q'}(\mathbb{R}_+;
L^{q'}(\mathbb{R}^n)))}.
\end{align}

\end{prop}

\begin{prop}\label{p2.9.1}
Let $2 \leq r \leq \infty$, $2/\nu(r) = n(1/2-1/r)$, $q > \nu>
\nu(r)\vee 2$, we have
\begin{align}
&\|\Box_k G_\nu u_0\|_{L^{\nu}_t(\mathbb{R}^+;
L^r_x(\mathbb{R}^n))}\lesssim \|\Box_k u_0\|_{L^2(\mathbb{R}^n)},\label{2.28}\\
&\|\Box_k \mathscr{A}_\nu f\|_{L^{\nu}_t(\mathbb{R}^+;
L^r_x(\mathbb{R}^n))\bigcap L^{\infty}_t(\mathbb{R}^+;
L^2_x(\mathbb{R}^n))}\lesssim C \|\Box_k f\|_{L^{q'}_t(\mathbb{R}^+;
L^{q'}_x(\mathbb{R}^n))}.\label{2.29}
\end{align}
\end{prop}
\noindent{\bf Proof.} Proposition \ref{prop3.3} implies \eqref{2.28}
and \eqref{2.29} directly.$\hfill\Box$

\vspace{10pt}

 Define the semigroup of Schr\"odinger equation
\begin{align}
 S(t)=\mathscr{F}^{-1}e^{-{\rm i}
 t\sum_{j=1}^{n}\xi_j^2}\mathscr{F}.\label{s3.3}
\end{align}

\begin{prop}\label{ubst}
$\Box_k G_\nu (t) : L^p \to L^p$ is uniformly bounded. More precisely,
\begin{align}
\|\Box_k G_\nu (t) u_0\|_{L^p(\mathbb{R}^n)}
 \lesssim (1+|t|^{n/2}) \|\Box_k  u_0\|_{L^p(\mathbb{R}^n)}
\label{modppa11}
\end{align}
uniformly holds for all $k\in \mathbb{Z}^n$, $\nu\ge 0$, $p\geq 1$.
\end{prop}

\noindent {\bf Proof.} It is well known that $e^{-|\xi|^2}$ is a
multiplier on $L^p$, i.e., $e^{-|\xi|^2}\in M_p$ ($M_p$ denotes
H\"ormander's multiplier space, see \cite{BJL}). Since $M_p$ is
isometrically invariant under affine transformations of
$\mathbb{R}^n$, we have $\|e^{-|\xi|^2}\|_{M_p} = \|e^{-\nu
t|\xi|^2}\|_{M_p}$, $1 \leq p \leq \infty$.
We have
\begin{align}
\|\Box_k G_\nu (t)f \|_p & \lesssim \|\Box_k S(t)f \|_p  \le \sum_{|\ell|_\infty \le 1} \|\mathscr{F}^{-1} \sigma_{k+\ell} e^{{\rm i}t |\xi|^2} \sigma_k \hat{f}\|_p \nonumber\\
& \le \sum_{|\ell|_\infty \le 1} \|\mathscr{F}^{-1} (\sigma_{k+\ell} e^{{\rm i}t |\xi|^2})\|_1 \| \Box_k f\|_p.
\label{modpp.1}
\end{align}
So, it suffices to show that $\|\mathscr{F}^{-1} (\sigma_{k} e^{{\rm i}t |\xi|^2})\|_1$ is uniformly bounded.
\begin{align}
\|\mathscr{F}^{-1} (\sigma_{k} e^{{\rm i}t |\xi|^2})\|_1
& = \|\mathscr{F}^{-1} (\sigma_{0} e^{{\rm i}t |\xi|^2})\|_1 \nonumber\\
& \lesssim \|\sigma_0\|^{1-n/2L}_2 \sum_{|\alpha|=L}\|D^\alpha (\sigma_{0} e^{{\rm i}t |\xi|^2})\|^{n/2L}_2 \nonumber\\
& \lesssim (1+|t|^{n/2}).
\label{modpp.2}
\end{align}
So we have the result, as desired. $\hfill \Box$

~\cite{WH} shows that $S(t)$ has the following maximal function
estimate:

\begin{lem}{\bf (\rm{\cite{WH}})} \label{3.9ss}
Let $4/n < p \leq \infty$, $p \geq 2$, $S(t)$ is defined as
\eqref{s3.3}, then we have
\begin{align}
\|\Box_k S(t)u_0\|_{L_{x_i}^{p}L_{x_j(j\neq
i)}^{\infty}L_t^{\infty}(\mathbb{R}\times \mathbb{R}^n)} \leq C
\langle k_i\rangle^{1/p} \|\Box_k u_0\|_{L^2(
\mathbb{R}^n)}.\label{3.9s}
\end{align}
\end{lem}
\begin{lem}\label{lem2.2}
Define maximal operator $M$ as following:
\begin{align}
(M f)(x) = \sup_{r>0}c_n r^{-n}\int_{|y|<r} |f(x-y)|d y.\nonumber
\end{align}
Let $\phi$ satisfies $\int_{\mathbb{R}^n} \phi \ dx = 1$, then for
any $f$, $f\in L_p$, $1<p\leq \infty$, we have
 \begin{align}
&\sup_{t>0}|f\ast \phi_t(x)|\leq M f(x) \int_{\mathbb{R}^n} \phi dx\label{2.9}\\
&\|M f\|_{L_p}\leq C \|f\|_{L_p}.\label{2.10}
 \end{align}
Where $\phi_t(x)= t^{-1}\phi(x/t)$.
\end{lem}
The proof can be found in  \cite{stein1},  Page 51,
\cite{stein2}, Page 3.

\begin{prop}\label{p2.5}
Let $4/n < p \leq \infty$, $p \geq 2$, we have
\begin{align}
\|\Box_k G_\nu(t)u_0\|_{L_{x_i}^{p}L_{x_j(j\neq
i)}^{\infty}L_t^{\infty}(\mathbb{R}^+\times \mathbb{R}^n)} \leq C
\langle k_i\rangle^{1/p} \|\Box_k u_0\|_{L^2(
\mathbb{R}^n)}.\label{3.9}
\end{align}
\end{prop}
\noindent{\bf Proof.} Take $i=1$ for example. When $t=0$,
\eqref{3.9} holds obviously. For $t>0$,
\begin{align}
\Box_k G_\nu(t)u_0&=\mathscr{F}^{-1}(e^{-\nu
t|\xi|^2})\ast\mathscr{F}^{-1}(e^{{-\rm i}t|\xi|^2}\widehat{\Box_k u_0})\nonumber\\
&=[\mathscr{F}^{-1}(e^{-|\xi|^2/2})]_{\sqrt{2\nu
t}}\ast\mathscr{F}^{-1}(e^{{-\rm i}t|\xi|^2}\widehat{\Box_k
u_0}).\label{2.9.1}
\end{align}
Notice that $\mathscr{F}^{-1}(e^{-|\xi|^2/2})=e^{-|x|^2/2}$,
$\int_{\mathbb{R}^n}e^{-|x|^2/2}d x = C$, then from \eqref{2.9},
\eqref{2.9.1} we have
\begin{align}
&\|\Box_k G_\nu(t)u_0\|_{L_{x_1}^{p}L_{x_j(j\neq
1)}^{\infty}(\mathbb{R}^n)L_{t\in(0,
\infty)}^{\infty}}\nonumber\\
&\leq \|[\mathscr{F}^{-1}(e^{-|\xi|^2/2})]_{\sqrt{2\nu
t}}\ast\mathscr{F}^{-1}(e^{{-\rm i}t'|\xi|^2}\widehat{\Box_k
u_0})\|_{L_{x_i}^{p}L_{x_j(j\neq
i)}^{\infty}(\mathbb{R}^n)L_{t'\in(0, \infty)}^{\infty}L_{t\in(0,
\infty)}^{\infty}}\nonumber\\
& \leq \|M [\mathscr{F}^{-1} (e^{{-\rm i}t|\xi|^2}\widehat{\Box_k
u_0})]\|_{L_{x_1}^{p}L_{x_j(j\neq
1)}^{\infty}(\times\mathbb{R}^n)L_{t\in(0, \infty)}^{\infty}}.
\end{align}
Define $M_{x_1}, M_{\bar{x}}$ were the maximal operators for
variable $x_1$ and the other varibles:
\begin{align}
&(M_{x_1} f)(x_1, \bar{x}) = \sup_{r>0}c_1 r^{-1}\int_{|y_1|<r}
|f(x_1-y_1, \bar{x})|d y_1,\nonumber\\
&(M_{\bar{x}} f)(x_1, \bar{x}) = \sup_{r>0}c_{n-1}
r^{-(n-1)}\int_{|\bar{y}|<r} |f(x_1, x_2-y_2,\ldots, x_n-y_n)|d
\bar{y}.\nonumber
\end{align}
From the definition of maximal operators and Lemma \ref{lem2.2} we
have
 \begin{align}
\|Mf(x_1, \bar{x})\|_{L_{x_1}^p L_{\bar{x}}^\infty}\leq \big\|
M_{x_1} \|M_{\bar{x}} f(x_1, \bar{x})\|_{
L_{\bar{x}}^\infty}\big\|_{L_{x_1}^p} &\leq \big\| M_{x_1} \| f(x_1,
\bar{x})\|_{
L_{\bar{x}}^\infty}\big\|_{L_{x_1}^p}\nonumber\\
&\leq \|f\|_{L_{x_1}^p L_{\bar{x}}^\infty},\label{4.37}
 \end{align}
where $p \geq 2$. From Lemma \ref{3.9ss} and \eqref{4.37} we obtain
 \begin{align}
  \|M  &  [\mathscr{F}^{-1} (e^{{-\rm i}t|\xi|^2}\mathscr{F}(\Box_k u_0))]\|
 _{L_{x_1}^{p}L_{x_j(j\neq
1)}^{\infty}L_t^{\infty}(\mathbb{R}^+\times
 \mathbb{R}^n)}\nonumber\\
&\leq C\|\mathscr{F}^{-1} (e^{{-\rm i}t|\xi|^2}\mathscr{F}(\Box_k
u_0))\|_{L_{x_1}^{p}L_{x_j(j\neq 1)}
 ^{\infty}L_t^{\infty}(\mathbb{R}\times \mathbb{R}^n)}\nonumber\\
 &\leq C \langle k_1\rangle^{1/p} \|\Box_k u_0\|_{L^2(
\mathbb{R}^n)}, \nonumber
\end{align}
which implies the result, as desired. $\hfill\Box$

\begin{prop}\label{p3.6}
For $n=1, 2$,  we have
\begin{align}
& \|\Box_k G_\nu(t)u_0\|_{L_{x}^{2} L_T^{\infty}} \leq C  \langle
k \rangle^{1/2} \ln 4\langle T\rangle  \|\Box_k
u_0\|_{L^1( \mathbb{R})}, \quad n=1; \nonumber\\
& \|\Box_k G_\nu(t)u_0\|_{L_{x_i}^{2}L_{x_j(j\neq
i)}^{\infty}L_T^{\infty}( \mathbb{R}^n)} \lesssim
       \langle k_i\rangle^{1/2}\|\Box_k u_0\|_{L^1( \mathbb{R}^n)}, \quad  n=2. \nonumber
\end{align}
\end{prop}
\noindent{\bf Proof.} We take $i=1$ for example.
\begin{align}
&\|\Box_k G_\nu(t)u_0\|_{L_{x_1}^{2}L_{x_j(j\neq 1)}^{\infty}L_T^{\infty}( \mathbb{R}^n)}\nonumber\\
&=\|\mathscr{F}_{\xi}^{-1}e^{-i t|\xi|^2}e^{-\nu
|t||\xi|^2}\tilde{\sigma}_k(\xi)\ast
\mathscr{F}_{\xi}^{-1}\sigma_k(\xi)
\widehat{u_0}\|_{L_{x_1}^{2}L_{x_j(j\neq 1)}^{\infty}L_T^{\infty}(
\mathbb{R}^n)}\nonumber\\
&\leq \|\mathscr{F}_{\xi}^{-1}e^{-i t|\xi|^2}e^{-\nu
|t||\xi|^2}\tilde{\sigma}_k(\xi)\|_{L_{x_1}^{2}L_{x_j(j\neq
1)}^{\infty}L_T^{\infty}( \mathbb{R}^n)} \|\Box_k
u_0\|_{L_{x_1}^{1}L_{x_j(j\neq 1)}^{1}( \mathbb{R}^n)},\label{3.16}
\end{align}
Where $\tilde{\sigma}_k(\xi)= \sum_{|l-k|< C(n)}\sigma_l(\xi)$. For
brevity, we still write $\tilde{\sigma}_k$ as $\sigma_k$. Now we
estimate $ \|\mathscr{F}_{\xi}^{-1}e^{-i t|\xi|^2}e^{-\nu
|t||\xi|^2}\sigma_k(\xi)\|_{L_{x_1}^{2}L_{x_j(j\neq
1)}^{\infty}L_T^{\infty}( \mathbb{R}^n)}$. First, consider the basic
$L^p-L^{p'}$ estimates for the semigroup of DGL equation $
G_{\nu}(t)= \mathscr{F}^{-1}e^{{-\rm i} t|\xi|^2}\cdot e^{-\nu
t|\xi|^2}\mathscr{F}$. We have
\begin{align}
\|G_\nu( t)\varphi\|_{L^\infty}\lesssim \|S( t)\varphi\|_{L^\infty}\lesssim  C
t^{-n/2}\|\varphi\|_{L^{1}},
\end{align}
and so,
\begin{align}
\|\mathscr{F}_\xi^{-1}e^{-{\rm i} t|\xi|^2}e^{-\nu
t|\xi|^2}\eta_{k_1}(\xi_1)\eta_{\bar{k}}(\bar{\xi})\|_{L_{x}^{\infty}
(\mathbb{R}^{n})}\leq C (1 + |t|)^{-n/2}.\label{g2.9}
\end{align}
On the other hand, using oscillatory integral techniques, we have
\begin{align}
\mathscr{F}_{\xi_1}^{-1}e^{-{\rm i} t\xi_1^2}e^{-\nu
t\xi_1^2}\eta_{k_1}(\xi_1)&=\int_\mathbb{R} e^{{\rm i}
x_1(\xi_1+\frac{t\xi_1^2}{x_1})}e^{-\nu t\xi_1^2}\eta_{k_1}(\xi_1)d
\xi_1\nonumber\\
&:=\int_\mathbb{R} e^{{\rm i} x_1 \phi(\xi_1)}\psi(\xi_1)d
\xi_1,\nonumber
\end{align}
where $\phi(\xi_1)= \xi_1 + \frac{t\xi_1^2}{x_1}$,
$\psi(\xi_1)=e^{-\nu t\xi_1^2}\eta_{k_1}(\xi_1)$. When
~$|x_1|>4\langle k_1\rangle|t|\vee 1$, we obtain $|\phi(\xi)'|>
1/2$. Meanwhile, it is easy to see
$$\int_{\mathbb{R}} |\psi(\xi_1)|d\xi_1 \leq C,\quad  \int_{\mathbb{R}}
|\psi'(\xi_1)| d\xi_1 \leq C, \quad  \int_{\mathbb{R}}
|\psi''(\xi_1)| d\xi_1 \leq C.$$ C is independent of $\nu, t$ and
$k_1$. So integrating by part, we obtain
\begin{align}
|\mathscr{F}_{\xi_1}^{-1}e^{-{\rm i} t\xi_1^2}e^{-\nu
t\xi_1^2}\eta_{k_1}(\xi_1)|\lesssim (1+ |x_1|)^{-2}.\label{g2.10.1}
\end{align}
from \eqref{g2.9}, \eqref{g2.10.1}, we have
\begin{align}
&\|\mathscr{F}^{-1}e^{-{\rm i}
t|\xi|^2}\eta_{k_1}(\xi_1)\eta_{\bar{k}}(\bar{\xi})\|
_{L_{x_1}^{2}L_{\bar{x}}^{\infty}L_{T}^{\infty}(\mathbb{R}^n)}\nonumber\\
&\leq \Big[\int_\mathbb{R}(1+ |x|)^{-4} dx_1\Big]^{1/2}+
\Big[\int_{|x_1|\leq 4\langle k_1\rangle T}\langle
k_1\rangle^{n}(\langle
k_1\rangle+|x_1|)^{-n}dx_1\Big]^{1/2}.
\end{align}
The result follows. $\hfill\Box$

Using similar method as in Proposition \ref{p3.6}, we
have
\begin{rem}\label{p3.6}
For $n\in \mathbb{N}$,  we have

$\|\Box_k G_\nu(t)u_0\|_{L_{x_i}^{1}L_{x_j(j\neq
i)}^{\infty}L_T^{\infty}( \mathbb{R}^n)} \leq$ C $\left\{
       \begin{array}{ll}
       \langle k_i\rangle \langle 4T \rangle^{1/2}\|\Box_k u_0\|_{L^1( \mathbb{R}^n)},                  & \hbox{$n=1$;} \\
       \langle k_i\rangle\ln \langle 4T \rangle \|\Box_k u_0\|_{L^1( \mathbb{R}^n)}, & \hbox{$n=2$;} \\
       \langle k_i\rangle\|\Box_k u_0\|_{L^1( \mathbb{R}^n)},                               & \hbox{$n\geq 3$.}
          \end{array}
        \right.$

\end{rem}

\vspace{10pt}

Next,  we  consider
the estimates between time-space norm and space-time norm
for integral operators $\mathscr{A}$. Since the semigroup of
Ginzburg-Landau equation does not have conjugate symmetry property as
Schr\"odinger equation, we can not apply $TT^*$ argument to obtain
some  good estimates as those of the Schr\"odinger equation, see
\cite{WH}.

\begin{prop}\label{p2.9}
Let $2\leq q < \infty$, $q> 4/n$, $\lambda=0, 1$, we have
\begin{align}
&\|\Box_k \mathscr{A}_\nu \partial_{x_i}^\lambda
f\|_{L_{x_i}^{2}L_{x_j(j\neq i)}^{\infty}L_t^{\infty}(\mathbb{R}^+
\times \mathbb{R}^n)} \lesssim \langle
k_i\rangle^{\lambda+1/2}\|\Box_k f\|_{L^{1}_t(\mathbb{R}^+;
L^{2}_x(\mathbb{R}^n))},\label{2.33}\\
&\|\Box_k \mathscr{A}_\nu \partial_{x_i}
f\|_{L_{x_i}^{\infty}L_{x_j(j\neq i)}^{2}L_t^{2}(\mathbb{R}^+ \times
\mathbb{R}^n)} \lesssim \langle k_i\rangle^{1/2}\|\Box_k
f\|_{L^{1}_t(\mathbb{R}^+; L^{2}_x(\mathbb{R}^n))}, \label{2.32}
\end{align}
where in \eqref{2.33},  condition $ |k_i|>4$ is required.
\end{prop}
\noindent{\bf Proof.} From \eqref{3.9}, \eqref{1.7}, Lemma
\ref{lem2.4} and Minkowski's inequality we have \eqref{2.33},
\eqref{2.32} hold, as desired. $\hfill \Box$

Similar to Proposition \ref{p2.9}, from Proposition\ref{p3.6}, we
have
\begin{prop}\label{p3.8}
For $n=1, 2$,  $T \leq 1$,  we have
\begin{align}
&\|\Box_k \mathscr{A}_\nu \partial_{x_i}
f\|_{L_{x_i}^{2}L_{x_j(j\neq i)}^{\infty}L_T^{\infty}(
\mathbb{R}^n)} \lesssim   \langle k_i\rangle^{3/2}   \|\Box_k
f\|_{L^{1}_T L^{1}_x(\mathbb{R}^n)}.\label{3.24}
\end{align}
\end{prop}

\begin{rem}\label{p3.10}
For $n\in \mathbb{N}$,  we have

$\|\Box_k \mathscr{A}_\nu \partial_{x_i}
f\|_{L_{x_i}^{1}L_{x_j(j\neq i)}^{\infty}L_T^{\infty}(
\mathbb{R}^n)}  \leq$ C $\left\{
       \begin{array}{ll}
       \langle k_i\rangle^{2} \langle 4T \rangle^{1/2} \|\Box_k f\|_{L^{1}_T L^{1}_x(\mathbb{R}^n)},                  & \hbox{n=1;} \\
       \langle k_i\rangle^{2}\ln \langle 4T \rangle \|\Box_k f\|_{L^{1}_T L^{1}_x(\mathbb{R}^n)}, & \hbox{n=2;} \\
       \langle k_i\rangle^{2}\|\Box_k f\|_{L^{1}_T L^{1}_x(\mathbb{R}^n)},                               & \hbox{n=3.}
          \end{array}
        \right.$

\end{rem}

From Propostions \ref{p2.9}, \ref{4.12} and \ref{prop2.1}, we
can obtain the following derivative interaction estimates:

\begin{lem}\label{3.1}
Let $i=2, \ldots, n$, we have
\begin{align}
&\|\Box_k \mathscr{A}_\nu\partial_{x_i}f\|_{L_{x_1}^{\infty}L_{x_2,
\ldots, x_n}^{2}L_t^2(\mathbb{R}^+ \times \mathbb{R}^n)} \le C
\|\partial_{x_i}\partial_{x_1}^{-1}\Box_k f\|_{L_{x_1}^{1}L_{x_2,
\ldots, x_n}^{2}L_t^2(\mathbb{R}^+
\times \mathbb{R}^n)},\label{3.1}\\
&\|\Box_k \mathscr{A}_\nu\partial_{x_i}f\|_{L_{x_1}^{\infty}L_{x_2,
\ldots, x_n}^{2}L_t^2(\mathbb{R}^+ \times \mathbb{R}^n)}\le C
\|\partial_{x_i}D_{x_1}^{-1/2}\Box_k f\|_{L^{1}_t(\mathbb{R}^+;
L^{2}_x(\mathbb{R}^n))},\label{3.2}\\
&\|\Box_k \mathscr{A}_\nu\partial_{x_i}f\|_{L_{x_1}^{2}L_{x_2,
\ldots, x_n}^{\infty}L_t^\infty(\mathbb{R}^+ \times \mathbb{R}^n)}
\le C\langle k_i\rangle\langle k_1\rangle^{1/2} \|\Box_k
f\|_{{L^{1}_t(\mathbb{R}^+; L^{2}_x(\mathbb{R}^n))}}.\label{3.4}
\end{align}
\end{lem}
Since the smooth-effect estimates for Ginzburg-Landau equation
\eqref{1.22} is almost the same with the Schr\"odinger equation (see
\cite{WH}). Follow the same method as \cite{WH}, we have

\begin{lem}\label{lem4.2}
Let $\psi: [0, \infty)\rightarrow [0,1]$ be a smooth bump function
satisfying $\psi(x)=1$ as $|x|\leq 1$ and $\psi(x)=0$ if $|x|\leq
2$. Denote $\psi_1(\xi)=\psi(\xi_2/2\xi_1)$,
$\psi_2(\xi)=1-\psi(\xi_2/2\xi_1)$, $\xi\in \mathbb{R}^n$. Then we
have for $\sigma \geq 0$,
\begin{align}
&\sum_{k\in\mathbb{Z}^n, |k_1|>4}\langle k_1\rangle^\sigma
\|\mathscr{F}_{\xi_1, \xi_2}^{-1}\psi_1 \mathscr{F}_{x_1, x_2}\Box_k
\partial_{x_2}\mathscr{A}_\nu f\|_{L_{x_1}^{\infty}L_{x_2,
\ldots, x_n}^{2}L_t^2(\mathbb{R}^+\times \mathbb{R}^n)}\nonumber\\
&\lesssim \sum_{k\in\mathbb{Z}^n, |k_1|>4}\langle k_1\rangle^\sigma
\|\Box_k f\|_{L_{x_1}^{1}L_{x_2, \ldots, x_n}^{2}L_t^2(\mathbb{R}_+
\times \mathbb{R}^n)},\label{4.5}
\end{align}
and for $\sigma\geq 1$,
\begin{align}
&\sum_{k\in\mathbb{Z}^n, |k_1|>4}\langle k_1\rangle^\sigma
\|\mathscr{F}_{\xi_1, \xi_2}^{-1}\psi_2 \mathscr{F}_{x_1, x_2}\Box_k
\partial_{x_2}\mathscr{A}_\nu f\|_{L_{x_1}^{\infty}L_{x_2,
\ldots, x_n}^{2}L_t^2(\mathbb{R}^+\times \mathbb{R}^n)}\nonumber\\
&\lesssim \sum_{k\in\mathbb{Z}^n, |k_2|>4}\langle k_2\rangle^\sigma
\|\Box_k f\|_{L_{x_1}^{1}L_{x_2, \ldots, x_n}^{2}L_t^2(\mathbb{R}_+
\times \mathbb{R}^n)}.\label{4.6}
\end{align}
\end{lem}

\section{Global well-posedness results for $n\geq 3$}

In this section, we will give the details of the proof of Theorem
\ref{thg1}. Define
\begin{align}
&\rho_1(u)=\sum_{i=1}^n  \sum_{k\in \mathbb{Z}^n, |k_i|>4}\langle
k_i\rangle^{s-1/2} \|\Box_k u\|_{L_{x_i}^{\infty}L_{(x_j)j\neq i
}^{2}L_t^2(\mathbb{R}^+\times \mathbb{R}^n)}:=\sum_{i=1}^n \rho_1^i(u) \nonumber\\
&\rho_2(u)=\sum_{i=1}^n \sum_{k\in \mathbb{Z}^n} \|\Box_k
u\|_{L_{x_i}^{2}L_{(x_j)j\neq i }^{\infty}L_t^\infty(\mathbb{R}_+
\times
\mathbb{R}^n)}\nonumber\\
&\rho_3(u)=\sum_{k\in \mathbb{Z}^n}\langle
k\rangle^{s-3/2  }\|\Box_k u\|_{L_t^\infty L_x^2\cap
L_t^3L_x^6(\mathbb{R}^+\times \mathbb{R}^n)}\nonumber
\end{align}
   Define resolution space as following:
\begin{align}
X_s:= \{u\in \mathscr{S}'(\mathbb{R}^+\times \mathbb{R}^n):
\|u\|_{X_s}:= \sum_{l=1}^3\sum_{\lambda=0,
1}\sum_{j=1}^n\rho_l(\partial_{x_j}^\lambda u)\leq
\delta_0\}\label{X}
\end{align}

\noindent{\bf{Proof of Theorem \ref{thg1}:}}  Using Lemma
\ref{appb}, we have for any $s>3$, there exist $\theta, \theta'>0$ such that
\begin{align}
\|u_0\|_{M_{2, 1}^3}\leq C \|u_0\|_{M_{2,
1}^{s }}^{1-\theta'}\|u_0\|_{L^2}^\theta.\label{4.2}
\end{align}
With the conditions that $u_0\in M_{2, 1}^{s }$, $\|u_0\|_{L^2}$
small enough, we can obtain $u_0 \in M_{2, 1}^{ 3}$ and
$\|u_0\|_{M_{2, 1}^{3}}$ sufficiently small.

 We only prove the result for the case $s=3$, we write $X_3=X$ for short.
Considering the following mapping:
\begin{align}
\mathscr{T}: u(t)\rightarrow G_\nu (t)u_0 -i
\mathscr{A}[\overrightarrow{\lambda_1}\cdot\nabla(|u|^2u) +
(\overrightarrow{\lambda_2}\cdot\nabla u)|u|^2+ \alpha
|u|^{2\delta}u], \nonumber
\end{align}
from \eqref{1.7}, \eqref{3.9}, \eqref{2.28} and Lemma \ref{lem2.4}
we have
\begin{align}
&\rho_1(\partial_{x_j}^\lambda G_\nu (t)u_0)\leq \sum_{i=1}^n
\sum_{k\in \mathbb{Z}^n, |k_i|>4}\langle k_i\rangle^{2}\langle
k_j\rangle^{\lambda} \|\Box_k u_0\|_{L_2(\mathbb{R}^n)}\lesssim
\|u_0\|_{M_{2, 1}^{3}}, \nonumber\\
&\rho_2(\partial_{x_j}^\lambda G_\nu (t)u_0)\leq \sum_{i=1}^n
\sum_{k\in \mathbb{Z}^n}\langle k_i\rangle^{1/2}\langle
k_j\rangle^{\lambda} \|\Box_k u_0\|_{L_2(\mathbb{R}^n)}\lesssim
\|u_0\|_{M_{2, 1}^{3/2}}, \nonumber\\
&\rho_3(\partial_{x_j}^\lambda G_\nu (t)u_0)\leq \sum_{k\in
\mathbb{Z}^n}\langle k \rangle^{3/2}\langle k_j\rangle^{\lambda}
\|\Box_k u_0\|_{L_2(\mathbb{R}^n)}\lesssim \|u_0\|_{M_{2,
1}^{5/2}}. \nonumber
\end{align}
So,  we obtain that
\begin{align}
\|G_\nu (t)u_0\|_X \lesssim \|u_0\|_{M_{2, 1}^{3}}. \label{5.1}
\end{align}
For the estimate of the nonlinear terms, noticing that
\begin{align}
\overrightarrow{\lambda_1}\cdot\nabla(|u|^2u) +
(\overrightarrow{\lambda_2}\cdot\nabla u)|u|^2 =\sum_{i=1}^n
[\lambda_1^i(\partial_{x_i}\bar{u}) u^2+ (2\lambda_1^i+\lambda_2^i)
(\partial_{x_i}u)u\bar{u}],\nonumber
\end{align}
and $\|u\|_X=\|\bar{u}\|_X$,  we only need to estimate
$\|\mathscr{A}_\nu(\lambda_1^i(\partial_{x_i}\bar{u}) u^2)\|_X$ and
$\|\mathscr{A}_\nu(\alpha |u|^{2\delta}u)\|_X$.
\begin{lem}{\rm (\cite{WaHe}, Lemma 7.2)}\label{lem5.1}
Let $s\geq 0, p \geq 1,\quad  p_i, \gamma, \gamma_i \leq \infty$
satisfy
\begin{align}
\frac{1}{p}=\frac{1}{p_1}+ \ldots + \frac{1}{p_N}, \quad
\frac{1}{\gamma}=\frac{1}{\gamma_1}+ \ldots +
\frac{1}{\gamma_N},\nonumber
\end{align}
then we have
\begin{align}
\sum_{k\in \mathbb{Z}^n}\langle k\rangle^{s}\|\Box_k(u_1 \ldots
u_N)\|_{L_t^\gamma L_x^p(\mathbb{R}^{1+n})} \lesssim \prod_{i=1}^N
\Big(\sum_{k\in \mathbb{Z}^n}\langle k\rangle^{s}\|\Box_k u_i
\|_{L_t^{\gamma_i} L_x^{p_i}(\mathbb{R}^{1+n})}\Big).
\end{align}
\end{lem}

From \eqref{2.33}, Lemma \ref{lem5.1} and H\"older's inequality we
have
\begin{align}
&\rho_2(\partial_{x_j}^\lambda(\mathscr{A}_\nu(\lambda_1^i(\partial_{x_i}\bar{u})
u^2))) \nonumber\\
&\lesssim \sum_{l=1}^n \sum_{k\in \mathbb{Z}^n} \|\Box_k
\partial_{x_j}^\lambda(\mathscr{A}_\nu(\partial_{x_i}\bar{u})
u^2)\|_{L_{x_l}^{2}L_{(x_m)m\neq l }^{\infty}L_t^\infty(\mathbb{R}_+
\times
\mathbb{R}^n)}\nonumber\\
&\lesssim n \sum_{k\in \mathbb{Z}^n} \langle
k_j\rangle^{\lambda}\langle k \rangle^{1/2}\|\Box_k
((\partial_{x_i}\bar{u}) u^2)\|_{L_{t}^1 L_x^2(\mathbb{R}^+\times
\mathbb{R}^n)}\nonumber\\
&\lesssim n \Big[ \sum_{k\in \mathbb{Z}^n} \langle
k\rangle^{3/2}\|\Box_k (\partial_{x_i}\bar{u})\|_{L_{t}^3
L_x^6(\mathbb{R}^+\times \mathbb{R}^n)}\times \big(\sum_{k\in
\mathbb{Z}^n} \langle k\rangle^{3/2}\|\Box_k u\|_{L_{t}^3
L_x^6(\mathbb{R}^+\times
\mathbb{R}^n)}\big)^2\Big]\nonumber\\
&\lesssim n \rho_3(\partial_{x_i}u)\rho_3(u)^2\label{5.4.1}
\end{align}
Noticing that $\|\Box_k u\|_{L_t^{k+1}L_x^{2(k+1)}(\mathbb{R}^+\times
\mathbb{R}^n)} \lesssim \|\Box_k u\|_{L_t^\infty L_x^2\cap
L_t^3L_x^6(\mathbb{R}^+\times \mathbb{R}^n)}$, follow the same
process as \eqref{5.4.1}, we have
\begin{align}
\rho_2(\partial_{x_j}^\lambda \mathscr{A}_\nu(\alpha
u^{\delta+1}\bar{u}^{\delta}))\lesssim \rho_3(u)^{2\delta
+1}\label{5.4.2}
\end{align}
From \eqref{5.4.1}, \eqref{5.4.2}, we obtain
\begin{align}
\rho_2\big[\partial_{x_j}^\lambda
\mathscr{A}_\nu(\lambda_1^i(\partial_{x_i}\bar{u}) u^2+\alpha
|u|^{2\delta}u)\big]\lesssim \|u\|_X^3 +
\|u\|_X^{2\delta+1}.\label{4.5.1}
\end{align}

Next, denote
\begin{align}
&\mathbb{S}_{\ell, 1}^{(i)} := \{(k^{(1)}, \ldots, k^{(\ell+1)})\in
(\mathbb{Z}^n)^{(\ell+1)}: |k_i^1|\vee \ldots \vee
|k_i^{(\ell+1)}|>4\},\nonumber\\
&\mathbb{S}_{\ell, 2}^{(i)} := \{(k^{(1)}, \ldots, k^{(\ell+1)})\in
(\mathbb{Z}^n)^{(\ell+1)}: |k_i^1|\vee \ldots \vee
|k_i^{(\ell+1)}|\leq 4\}.\nonumber
\end{align}
Using the frequency-uniform decomposition, we have
\begin{align}
\prod_{r=1}^{\ell+1} u_r
= \sum_{\mathbb{S}_{\ell, 1}^{(i)}} \Box_{k^{(1)}}u_1 \ldots
\Box_{k^{(\ell+1)}}u_{\ell+1} + \sum_{\mathbb{S}_{\ell, 2}^{(i)}}
\Box_{k^{(1)}}u_1 \ldots \Box_{k^{(\ell+1)}}u_{\ell+1}.\label{5.2}
\end{align}
Where we divide $(\mathbb{Z}^n)^{(\ell+1)}$ into two parts
$\mathbb{S}_{ \ell,1}^{(i)}$ and $\mathbb{S}_{\ell, 2}^{(i)}$ by
considering variable $x_i$. $\mathbb{S}_{\ell, 1}^{(i)}$ denotes the
high frequency part, so we will apply smooth effect estimates; while
$\mathbb{S}_{\ell, 2}^{(i)}$ denotes the low frequency part, we will
apply Strichartz-type estimates.

Now we estimate
$\rho_1(\partial_{x_j}^\lambda(\mathscr{A}_\nu\lambda_1^i(\partial_{x_i}\bar{u})
u^2)$. It suffices to show that
\begin{align}
\rho_1(\partial_{x_j}^\lambda\mathscr{A}_\nu
(\lambda_1^i(\partial_{x_i}\bar{u}) u^2))\lesssim \|u\|_X^{3}, \quad
j= 1, 2.\nonumber
\end{align}
 The other cases are almost the same.
For the estimates of
$\rho_1^1(\partial_{x_j}^\lambda(\mathscr{A}_\nu\lambda_1^i(\partial_{x_i}\bar{u})
u^2)) $ are similar to the proof in \cite{WH}, we leave the details
of the proof into Appendix A.

collecting \eqref{5.3}--\eqref{5.10}, from symmetry, we can obtain
that
\begin{align}
\sum_{\lambda=0,
1}\sum_{j=1}^n\rho_1(\partial_{x_j}^\lambda\mathscr{A}_\nu
(\lambda_1^i(\partial_{x_i}\bar{u}) u^2))\leq
\|u\|_X^{3}.\label{5.11}
\end{align}

The estimate of $\rho_1(\partial_{x_j}^\lambda\mathscr{A}_\nu(\alpha
|u|^{2\delta}u))$ is similar to
$\rho_1(\partial_{x_j}^\lambda\mathscr{A}_\nu((\partial_{x_i}\bar{u})
u^2))$, but in \eqref{5.5}, we will apply
\begin{align}
&\|\Box_{k^{(1)}}u\ldots\Box_{k^{(2\sigma+1)}}u\|_{L_{x_1}^{1}L_{x_2,
\ldots, x_n}^{2}L_t^2(\mathbb{R}^+ \times \mathbb{R}^n)}
\nonumber\\
& \lesssim \|\Box_{k^{(1)}} u\|_{L_{x_1}^{\infty}L_{x_2, \ldots,
x_n}^{2}L_t^2(\mathbb{R}^+ \times
\mathbb{R}^n)}\prod_{i=2}^{2\sigma+1}\|\Box_{k^{(i)}}
u\|_{L_{x_1}^{2}L_{x_2, \ldots, x_n}^{\infty}L_t^\infty\cap
L_t^\infty L_x^2(\mathbb{R}^+ \times \mathbb{R}^n)}.\label{4.70}
\end{align}
In \eqref{5.6}, we will apply
\begin{align}
&\|\Box_k u\|_{L_t^pL_x^{2p}(\mathbb{R}^+ \times \mathbb{R}^n)}
 \lesssim\|\Box_k u\|_{L_t^\infty L_x^2\cap
L_t^3L_x^6(\mathbb{R}^+ \times \mathbb{R}^n)}, \quad p\geq
3.\label{4.71}
\end{align}

Finally, we consider the estimate of
$\rho_3(\partial_{x_j}^\lambda\mathscr{A}_\nu(
\lambda_1^i(\partial_{x_i}\bar{u}) u^2))$. From Lemma \ref{lem2.4}
and \eqref{2.29} (where we let $\nu=3, r=6, q=2\sigma + 2$, $\sigma
\geq 1$). we can obtain
\begin{align}
\|\Box_k \partial_{x_j}^\lambda\mathscr{A}_\nu
(\lambda_1^i(\partial_{x_i}\bar{u}) u^2)\|_{L_t^\infty L_x^2\cap
L_t^3L_x^6(\mathbb{R}^+ \times \mathbb{R}^n)}&\lesssim
\|\Box_k(\partial_{x_j}^\lambda (\lambda_1^i(\partial_{x_i}\bar{u})
u^2))\|_{L_{t, x}^{\frac{2\sigma+2}{2\sigma+1}}(\mathbb{R}^+ \times \mathbb{R}^n)}\nonumber\\
&\lesssim \langle k_j\rangle^\lambda \|\Box_k
((\partial_{x_i}\bar{u})u^2)\|_{L_{t,
x}^{\frac{2\sigma+2}{2\sigma+1}}(\mathbb{R}^+ \times
\mathbb{R}^n)}.\label{5.15}
\end{align}
 Let $\sigma=1$ in \eqref{5.15}, applying Young's inequality we have:
\begin{align}
&\rho_3(\partial_{x_j}^\lambda\mathscr{A}_\nu
(\lambda_1^i(\partial_{x_i}\bar{u}) u^2))\nonumber\\
&\lesssim \sum_{k\in \mathbb{Z}^n}\langle
 k\rangle^{3/2}\langle
 k_j\rangle^{\lambda}
\sum_{k^{(1)},k^{(2)}, k^{(3)}}\|\Box_k(\Box_{k^{(1)}}
(\partial_{x_i}\bar{u}) \Box_{k^{(2)}}u\Box_{k^{(3)}}u)\|_{L_{t,
x}^{4/3}(\mathbb{R}^+ \times \mathbb{R}^n)}\nonumber \\
&\lesssim \sum_{k\in \mathbb{Z}^n}(\langle
 k\rangle^{5/2}+\langle
 k_j\rangle^{5/2})
\sum_{k^{(1)},k^{(2)}, k^{(3)}}\|\Box_k(\Box_{k^{(1)}}
(\partial_{x_i}\bar{u}) \Box_{k^{(2)}}u\Box_{k^{(3)}}u)\|_{L_{t,
x}^{4/3}(\mathbb{R}^+ \times \mathbb{R}^n)}\nonumber\\
&\lesssim \sum_{j=1}^n\sum_{k\in \mathbb{Z}^n}\langle
 k_j\rangle^{5/2}
\sum_{k^{(1)},k^{(2)}, k^{(3)}}\|\Box_k(\Box_{k^{(1)}}
(\partial_{x_i}\bar{u}) \Box_{k^{(2)}}u\Box_{k^{(3)}}u)\|_{L_{t,
x}^{4/3}(\mathbb{R}^+ \times \mathbb{R}^n)}\nonumber
\end{align}
Then applying the decomposition \eqref{5.2}(where we consider the
variable $x_j$),
 and obtain
\begin{align}
&\sum_{k\in \mathbb{Z}^n}\langle
 k_j\rangle^{5/2}
\sum_{k^{(1)},k^{(2)}, k^{(3)}}\|\Box_k(\Box_{k^{(1)}}
(\partial_{x_i}\bar{u}) \Box_{k^{(2)}}u\Box_{k^{(3)}}u)\|_{L_{t,
x}^{4/3}(\mathbb{R}^+ \times \mathbb{R}^n)}\nonumber\\
&\lesssim
 \sum_{k\in \mathbb{Z}^n, |k_j|>4}
\langle
 k_j\rangle^{5/2}
\sum_{\mathbb{S}_{2, 1}^{(j)}}\|\Box_k(\Box_{k^{(1)}}
(\partial_{x_i}\bar{u}) \Box_{k^{(2)}}u\Box_{k^{(3)}}u)\|_{L_{t,
x}^{4/3}(\mathbb{R}^+ \times \mathbb{R}^n)}\nonumber\\
& + \sum_{k\in \mathbb{Z}^n, |k_j|>4}\langle
 k_j\rangle^{5/2}
\sum_{\mathbb{S}_{2, 2}^{(j)}}\|\Box_k(\Box_{k^{(1)}}
(\partial_{x_i}\bar{u}) \Box_{k^{(2)}}u\Box_{k^{(3)}}u)\|_{L_{t,
x}^{4/3}(\mathbb{R}^+ \times
\mathbb{R}^n)}\nonumber\\
& + \sum_{k\in \mathbb{Z}^n, |k_j|\leq 4}\langle
 k_j\rangle^{5/2}
\sum_{k^{(1)},k^{(2)}, k^{(3)}\in
\mathbb{Z}^3}\|\Box_k(\Box_{k^{(1)}} (\partial_{x_i}\bar{u})
\Box_{k^{(2)}}u\Box_{k^{(3)}}u)\|_{L_{t,
x}^{4/3}(\mathbb{R}^+ \times \mathbb{R}^n)}\nonumber \\
& = V + VI + VII.\label{5.16}
\end{align}

For the term $VI, VII$, from definition of $\mathbb{S}_{2, 2}^{(j)}$
, we can see $\langle
 k_j\rangle\leq C$.
Follow the same process of \eqref{5.6}, we have
\begin{align}
&VI + VII \nonumber\\
 &\lesssim
\sum_{k^{(1)},k^{(2)}, k^{(3)}\in \mathbb{Z}^n}\|\Box_{k^{(1)}}
(\partial_{x_i}\bar{u}) \Box_{k^{(2)}}u\Box_{k^{(3)}}u\|_{L_{t,
x}^{4/3}(\mathbb{R}^+ \times \mathbb{R}^n)}\nonumber \\
&\lesssim \sum_{k^{(1)},k^{(2)}, k^{(3)}\in \mathbb{Z}^n}
\|\Box_{k^{(1)}}
\partial_{x_i}\bar{u}
\|_{ L_t^4L_x^4(\mathbb{R}^+ \times\mathbb{R}^n)}\|\Box_{k^{(2)}}u
\|_{ L_t^4L_x^4(\mathbb{R}^+ \times\mathbb{R}^n)}\|\Box_{k^{(3)}}u
\|_{ L_t^4L_x^4(\mathbb{R}^+ \times\mathbb{R}^n)}\nonumber\\
&\lesssim \sum_{k^{(1)}\in \mathbb{Z}^n}\langle
k^{(1)}\rangle^{3/2}\|\Box_{k^{(1)}}
\partial_{x_i}\bar{u}
\|_{ L_t^4L_x^4(\mathbb{R}^+ \times\mathbb{R}^n)}(\sum_{k^{(2)}\in
\mathbb{Z}^n}\langle k^{(2)}\rangle^{3/2}\|\Box_{k^{(2)}} u
\|_{ L_t^4L_x^4(\mathbb{R}^+ \times\mathbb{R}^n)})^2\nonumber\\
&\lesssim \rho_3(\partial_{x_i}u)\rho_3(u)^2.\label{4.74}
\end{align}
At the last step of \eqref{4.74}, we apply
\begin{align}
&\|\Box_k u\|_{L_t^pL_x^{p}(\mathbb{R}^+ \times \mathbb{R}^n)}
 \lesssim\|\Box_k u\|_{L_t^\infty L_x^2\cap
L_t^3L_x^6(\mathbb{R}^+ \times \mathbb{R}^n)}, \quad p\geq
4.\label{4.7a}
\end{align}
Now, we consider the estimate of $V$, from H\"older's inequality and
\eqref{4.71}, we can obtain
\begin{align}
&\|\Box_{k^{(1)}} (\partial_{x_i}\bar{u})
\Box_{k^{(2)}}u\Box_{k^{(3)}}u\|_{L_{t, x}^{4/3}(\mathbb{R}^+ \times
\mathbb{R}^n)}\nonumber\\
&\leq \|\Box_{k^{(1)}} (\partial_{x_i}\bar{u})
|\Box_{k^{(2)}}u\Box_{k^{(3)}}u|^{1/2}\|_{L_{t, x}^{2}(\mathbb{R}^+
\times \mathbb{R}^n)}\|
|\Box_{k^{(2)}}u\Box_{k^{(3)}}u|^{1/2}\|_{L_{t, x}^{4}(\mathbb{R}^+
\times
\mathbb{R}^n)}\nonumber\\
&\leq\|\Box_{k^{(1)}} (\partial_{x_i}\bar{u})\|
_{L_{x_j}^{\infty}L_{x_r(r\neq j)}^{2}L_t^2(\mathbb{R}^+ \times
\mathbb{R}^n)} \|\Box_{k^{(2)}}u \|^{1/2}_{L_{x_j}^{2}L_{x_r(r\neq
j)}^{\infty}L_t^\infty(\mathbb{R}^+ \times \mathbb{R}^n)} \|
\Box_{k^{(3)}}u \|^{1/2}_{L_{x_j}^{2}L_{x_r(r\neq
j)}^{\infty}L_t^\infty(\mathbb{R}^+ \times \mathbb{R}^n)}\nonumber\\
&\quad\quad\times \|\Box_{k^{(2)}} u\|^{1/2}_{L_t^\infty L_x^2\cap
L_t^3 L_x^6(\mathbb{R}^+ \times\mathbb{R}^n)} \| \Box_{k^{(3)}}u
\|^{1/2}_{L_t^\infty L_x^2\cap L_t^3 L_x^6(\mathbb{R}^+
\times\mathbb{R}^n)}.\label{4.77}
\end{align}
In this way the estimate of $V$ reduces to the estimate of
$I$ as in \eqref{5.5}. Since $|k_j-k_j^{(1)}-k_j^{(2)}-k_j^{(3)}| \leq
C$, without loss of generality, we can assume $
 |k_j^{(1)}|=\max_{r=1,2,3}|k_j^r|$, from \eqref{4.77} and H\"older's inequality we have
\begin{align}
V
&\lesssim \sum_{k^{(1)}\in \mathbb{Z}^n, |k_j^{(1)}|>4}\langle
k_j^{(1)}\rangle^{5/2} \|\Box_{k^{(1)}}
\partial_{x_i}\bar{u}\|_{L_{x_j}^{\infty}L_{x_r(r\neq j)}^{2}L_t^2(\mathbb{R}^+
\times \mathbb{R}^n)}\nonumber\\
&\quad \quad\times \sum_{k^{(2)}\in \mathbb{Z}^n}
\Big(\|\Box_{k^{(2)}} u\|_{L_{x_j}^{2}L_{x_r(r\neq
j)}^{\infty}L_t^\infty(\mathbb{R}^+ \times \mathbb{R}^n)}+
\|\Box_{k^{(2)}} u\|_{L_t^\infty L_x^2\cap L_t^3
L_x^6(\mathbb{R}^+ \times\mathbb{R}^n)}\Big)\nonumber\\
&\quad\quad \times \sum_{k^{(3)}\in\mathbb{Z}^n}
 \Big(\|\Box_{k^{(3)}}u\|_{L_{x_j}^{2}L_{x_r(r\neq j)}^{\infty}L_t^\infty(\mathbb{R}^+ \times
\mathbb{R}^n)}+ \|\Box_{k^{(3)}}u\|_{L_t^\infty L_x^2\cap
L_t^3L_x^6(\mathbb{R}^+\times\mathbb{R}^n)}\Big)\nonumber\\
&\lesssim \rho_1^j({\partial_{x_i}u})(\rho_2(u) +
\rho_3(u))^2.\label{4.78}
\end{align}
The estimate for $\rho_3(\partial_{x_j}^\lambda\mathscr{A}(\alpha
|u|^{2\delta}u))$ is similar to
$\rho_3(\partial_{x_j}^\lambda\mathscr{A}((\partial_{x_i}\bar{u})
u^2))$, the difference is that in \eqref{4.77}, we will use
\begin{align}
&\|\Box_{k^{(1)}}u\ldots\Box_{k^{(\sigma+1)}}u\Box_{k^{(\sigma+2)}}\bar{u}\ldots\Box_{k^{(2\sigma+1)}}\bar{u}\|_{L_{t,
x}^{\frac{2\sigma+2}{2\sigma+1}}(\mathbb{R}^+ \times
\mathbb{R}^n)}\nonumber\\
&\leq\|\Box_{k^{(1)}}u\|_{L_{x_j}^{\infty}L_{x_r(r\neq
j)}^{2}L_t^2(\mathbb{R}^+ \times \mathbb{R}^n)}
\prod_{m=2}^{\sigma+1}\|\Box_{k^{(m)}}u
\|^{\sigma}_{L_{x_j}^{2\sigma}L_{x_r(r\neq
j)}^{\infty}L_t^\infty(\mathbb{R}^+ \times \mathbb{R}^n)} \nonumber\\
& \quad \times
\prod_{m=\sigma+2}^{2\sigma+1}\|\Box_{k^{(m)}}
\bar{u}\|^{\sigma}_{L_{t, x}^{2\sigma+1}(\mathbb{R}^+
\times\mathbb{R}^n)}.\nonumber
\end{align}
In addition to \eqref{4.7a} and \eqref{4.70}, we can obtain the
estimate of $\rho_3(\partial_{x_j}^\lambda\mathscr{A}(\alpha
|u|^{2\delta}u))$. Until now, we have obtain
\begin{align}
\rho_3(\partial_{x_j}^\lambda\mathscr{A}_\nu(\alpha
|u|^{2\delta}u)+\rho_3(\partial_{x_j}^\lambda\mathscr{A}_\nu
(\lambda_1^i(\partial_{x_i}\bar{u}) u^2)))\lesssim
\|u\|_X^{2\delta+1}+\|u\|_X^{3}.\label{4.73}
\end{align}

Collecting \eqref{4.5.1}, \eqref{5.11}, \eqref{4.73}
\begin{align}
\|\mathscr{A}_\nu u\|_X \lesssim \|u_0\|_{M_{2, 1}^{3}} + n
(\|u\|_X^{3}+\|u\|_X^{2\delta+1}).
\end{align}
Using standard contraction mapping argument, we can obtain that Eq. \eqref{GL} has a unique solution $u \in X$ with $\|u\|_X \le C \|u_0\|_{M^3_{2,1}}$.

Finally, for the general case $s>3$,  using similar way as in the above, we have
\begin{align}
\|u\|_{X_{s }} \lesssim \|u_0\|_{M_{2,
1}^{s }}+\|u\|_{X_{s }}\|u\|_{X_{3}}^2
+\|u\|_{X_{s }}\|u\|_{X_{3}}^{2\delta}.\label{4.20}
\end{align}
Since in the right hand side of \eqref{4.20}, Using the fact that $\|u\|_{X_{3}}$ is sufficiently small, we can get that
\begin{align}
\|u\|_{X_{s }} \lesssim \|u_0\|_{M_{2,
1}^{s }}.\label{4.21}
\end{align}
Finally, we show that $u\in L^\infty_t (\mathbb{R}^{+}; M_{2, 1}^s
(\mathbb{R}^n))$. From Proposition \ref{4.12} and Proposition
\ref{p2.9}, we have
\begin{align}
& \sum_{k\in \mathbb{Z}^n}\langle
 k\rangle^{s}\|\Box_k \mathscr{A}_\nu[\overrightarrow{\lambda_1}\cdot\nabla(|u|^2u) +
(\overrightarrow{\lambda_2}\cdot\nabla
u)|u|^2]\|_{L^\infty_tL_x^2(\mathbb{R}^+ \times
 \mathbb{R}^n)}\nonumber\\
& \lesssim  \sum_{j=1}^n
 \sum_{k\in \mathbb{Z}^n, |k_j|>4}\langle
 k_j\rangle^{s-1/2}
\sum_{\mathbb{S}^{(1)}_{2,1}}\|\Box_k (\Box_{k^{(1)}}
\partial_{x_i}\bar{u}
\Box_{k^{(2)}}u\Box_{k^{(3)}}u)\|_{L_{x_j}^{1}L_{x_i(i\neq
j)}^{2}L_t^2(\mathbb{R}_+
\times \mathbb{R}^n)}\nonumber\\
& \quad + \sum_{j=1}^n
 \sum_{k\in \mathbb{Z}^n, |k_j|>4}
\sum_{\mathbb{S}^{(1)}_{2,2}}\|\Box_k (\Box_{k^{(1)}}
\partial_{x_i}\bar{u}
\Box_{k^{(2)}}u\Box_{k^{(3)}}u)\|_{L_{t}^{1}L_{x}^{2} (\mathbb{R}_+
\times \mathbb{R}^n)}\nonumber\\
&\quad + \sum_{j=1}^n
 \sum_{k\in \mathbb{Z}^n, |k_j|\leq 4}\langle
 k_j\rangle^{s}
\sum_{k^{(1)}, k^{(2)}, k^{(3)}}\|\Box_k (\Box_{k^{(1)}}
\partial_{x_i}\bar{u}
\Box_{k^{(2)}}u\Box_{k^{(3)}}u)\|_{L_{t}^{1}L_{x}^{2} (\mathbb{R}_+
\times \mathbb{R}^n)}.\label{4.23}
\end{align}
From \eqref{4.23} and the estimate of $\rho_1$ and part $IIV$ in
$\rho_3$ above, we obtain
\begin{align}
 \sum_{k\in \mathbb{Z}^n}\langle
 k\rangle^{s}\|\Box_k u\|_{L^\infty_tL_x^2(\mathbb{R}^+ \times
 \mathbb{R}^n)}
 \lesssim \|u_0\|_{M_{2, 1}^s} +  \|u\|_{X_s}^3+
\|u\|_{X_s}^{2\delta+1},
\end{align}
which implies that $u\in L^\infty_t (\mathbb{R}^{+}; M_{2, 1}^s
(\mathbb{R}^n))$.

\section {Limit behavior as $\nu \rightarrow 0$}

 In this section, we will prove Theorem \ref{thg2}. Letting
initial data $u_0$ belong to $ M_{2,1}^{4}$ and $\|u_0\|_{L^2}$
small enough, we prove that the solution of derivative
Ginzburg-Landau equation \eqref{GL} will converge to that of
derivative Schr\"odinger equation \eqref{nls}
 as $\nu \rightarrow 0$.

 Let
$S(t)=\mathscr{F}^{-1}e^{{-\rm i} t|\xi|^2}
\mathscr{F}$ denote the semi-group of derivative Schr\"odinger
equation and $\mathscr{L}f(t, x)=\int_0^t S(t-\tau)f(\tau, x)d
\tau$. Rewrite DCGL equation \eqref{GL} as
\begin{align}
u= S(t)u_0 +
\mathscr{L}[\overrightarrow{\lambda_1}\cdot\nabla(|u|^2u) +
(\overrightarrow{\lambda_2}\cdot\nabla u)|u|^2+ \alpha
|u|^{2\delta}u] + \nu\mathscr{L}(\triangle u).
\end{align}
Then define:
\begin{align}
&\rho_T^1(u)=\sum_{i=1}^n  \sum_{k\in \mathbb{Z}^n, |k_i|>4}\langle
k_i\rangle  \|\Box_k u\|_{L_{x_i}^{\infty}L_{(x_j)j\neq i
}^{2}L_t^2([0, T] \times \mathbb{R}^n)}, \nonumber\\
&\rho_T^2(u)=\sum_{i=1}^n \sum_{k\in \mathbb{Z}^n} \|\Box_k
u\|_{L_{x_i}^{2}L_{(x_j)j\neq i }^{\infty}L_t^\infty([0, T] \times
\mathbb{R}^n)},\nonumber\\
&\rho_T^3(u)=\sum_{k\in \mathbb{Z}^n} \|\Box_k
u\|_{L_t^\infty L_x^2\cap L_t^3L_x^6([0, T] \times
\mathbb{R}^n)}.\nonumber
\end{align}
\begin{align}
\|u\|_{Y_T}:=\sum_{\lambda=0,
1}\sum_{j=1}^n(\rho_T^1(\partial_{x_j}^\lambda
u)+\rho_T^2(\partial_{x_j}^\lambda
u)+\rho_T^3(\partial_{x_j}^\lambda u)) \label{gy}
\end{align}
×¢Òâµ½
$$
\lambda_1\cdot\nabla(|u|^2u) + (\lambda_2\cdot\nabla u)|u|^2
=\sum_{i=1}^n [\lambda_1^i(\partial_{x_i}\bar{u}) u^2+
(2\lambda_1^i+\lambda_2^i) (\partial_{x_i}u)u\bar{u}].\nonumber
$$
Denote $v$ is the solution of derivative nonlinear Schr\"odinger
equation \eqref{nls} with the same initial data. Combining the method
in \cite{FDF}, we only need to estimate the following
\begin{align}
&\|u_\nu - v\|_{Y_T} &\nonumber\\
&\lesssim
\sum_{i=1}^n\|\mathscr{L}\lambda_1^i((\partial_{x_i}\bar{u}_\nu)
u_\nu^2 -(\partial_{x_i}\Bar{v})v^2)\|_{Y_T} +
\sum_{i=1}^n\|\mathscr{L}(2\lambda_1^i+\lambda_2^i)((\partial_{x_i}u_\nu)
|u_\nu|^2-(\partial_{x_i}v)|v|^2)\|_{Y_T}\nonumber\\
& \quad +\|\mathscr{L}(\bar{u}_\nu^\sigma u_\nu^{\sigma +1}-\bar{v}^\sigma
v^{\sigma +1})\|_{Y_T}+ \nu\|\mathscr{L}\triangle
u_\nu\|_{Y_T}.\label{g6.1}
\end{align}
Similar to the method in \cite{FDF}, take the first and third term
in \eqref{g6.1} for example. The second term can be treated in
similar way.
\begin{align}
&(\partial_{x_i}\bar{u}_\nu) u_\nu^2 -(\partial_{x_i}\Bar{v})v^2
=\partial_{x_i}(\bar{u}_\nu -\bar{v})u_\nu^2 +
\partial_{x_i}\bar{v}(u_\nu-v)(u_\nu+v),\label{4.88}\\
&\bar{u}_\nu^\sigma u_\nu^{\sigma +1}-\bar{v}^\sigma v^{\sigma +1}=
\bar{u}_\nu^\sigma(u_\nu-v)\sum_{q=0}^{\sigma}u_\nu^q v^{\sigma-q} +
{v}^{\sigma+1}(\bar{u}_\nu-\bar{v})\sum_{q=0}^{\sigma-1}\bar{u}_\nu^q
\bar{v}^{\sigma-1-q}.
\end{align}
Using the decomposition in \eqref{5.2} and combine the proof in
Section 4, we only need to substitute $\partial_{x_i}\bar{u}$ with
$\partial_{x_i}(\bar{u}_\nu -\bar{v})$ in the proof of
$\|\mathscr{L}(\lambda_1^i(\partial_{x_i}\bar{u}) u^2)\|_X$, then we
have
\begin{align}
\|\mathscr{L}\lambda_1^i(\partial_{x_i}(\bar{u}_\nu -\bar{v})
u_\nu^2) \|_{Y_T} \leq C \|u_\nu-v\|_{Y_T} \|u_\nu\|_X^2.\nonumber
\end{align}
Then substitute $u^2$ with $(u_\nu-v)(u_\nu+v)$, we have
\begin{align}
\|\mathscr{L}\lambda_1^i(\partial_{x_i}\bar{v}(u_\nu-v)(u_\nu+v))
\|_{Y_T} \leq C \|u_\nu-v\|_{Y_T} \|v\|_X(\|v\|_X+
\|u_\nu\|_X).\nonumber
\end{align}
Repeat the argument in \eqref{4.2}, we have $\|u_\nu\|_X, \|v\|_X$
are sufficiently small. Then we have
\begin{align}
\|\mathscr{L}\lambda_1^i((\partial_{x_i}\bar{u}_\nu) u_\nu^2
-(\partial_{x_i}\Bar{v})v^2)\|_{Y_T}\lesssim \frac{1}{10} \|u_\nu-v\|_{Y_T}
\end{align}
Similarly, in the estimate of $\|\alpha u^{\delta +1}u^\delta\|_Y$,
we will substitute $u^{\delta +1}u^\delta$ with
$\bar{u}_\nu^\delta(u_\nu-v)u_\nu^q v^{\delta-q}$ and
${v}^{\delta+1}(\bar{u}_\nu-\bar{v})\bar{u}_\nu^q
\bar{v}^{\delta-1-q}$, we have
\begin{align}
\|\mathscr{L}(\bar{u}_\nu^\delta u_\nu^{\delta +1}-\bar{v}^\delta
v^{\delta +1})\|_{Y_T} & \leq C \|u_\nu-v\|_{Y_T}
\big[\|u_\nu\|_X^{\delta+q}\|v\|_X^{\delta-q}+\|u_\nu\|_X^{q}\|v\|_X^{2\delta-q}\big] \nonumber\\
& \lesssim
\frac{1}{10} \|u_\nu-u\|_{Y_T}.\nonumber
\end{align}
 Moving the first three term in the right of
\eqref{g6.1} to the left, then from the definition of $Y_T$ we
obtain
\begin{align}
\|u_\nu-v\|_{Y_T}&\leq C \nu\|\mathscr{L}\triangle u\|_{Y_T}\nonumber\\
&\lesssim \sum_{\lambda=0, 1}\sum_{j=1}^n\sum_{i=1}^n  \sum_{k\in
\mathbb{Z}^n, |k_i|>4}\langle k_i\rangle
\|\partial_{x_j}^\lambda \Box_k (\mathscr{L}\triangle
u)\|_{L_{x_i}^{\infty}L_{(x_j)j\neq i
}^{2}L_t^2([0, T]\times \mathbb{R}^n)}\nonumber\\
&+\sum_{\lambda=0, 1}\sum_{j=1}^n\sum_{i=1}^n \sum_{k\in
\mathbb{Z}^n} \|\partial_{x_j}^\lambda\Box_k (\mathscr{L}\triangle
u)\|_{L_{x_i}^{2}L_{(x_j)j\neq i
}^{\infty}L_t^\infty([0, T] \times \mathbb{R}^n)}\nonumber\\
&+ \sum_{\lambda=0, 1}\sum_{j=1}^n\sum_{k\in \mathbb{Z}^n} \|\partial_{x_j}^\lambda\Box_k (\mathscr{L}\triangle
u)\|_{L_t^\infty L_x^2\cap
L_t^3L_x^6([0, T] \times \mathbb{R}^n)}\nonumber\\
& := A_1 + A_2 + A_3.\label{4.93}
\end{align}
Similar to \cite{FDF}, applying  Minkowski's inequality and
\eqref{2.28}, \eqref{3.9} and \eqref{1.7}, we can obtain
\begin{align}
&A_1\lesssim \nu T  \sum_{\lambda=0, 1}\sum_{j=1}^n\sum_{i=1}^n
\sum_{k\in \mathbb{Z}^n}  \langle k_i\rangle^{1/2}\langle
k\rangle^{2}\|\Box_k (\partial_{x_j}^\lambda u)\|_{L_T^\infty
L_x^2(\mathbb{R}^n)},\label{4.92}
\\
&A_2\lesssim \nu T  n \sum_{\lambda=0, 1}\sum_{j=1}^n \sum_{k\in
\mathbb{Z}^n} \langle k\rangle^{2}\langle k\rangle^{1/2}\|\Box_k
(\partial_{x_j}^\lambda u)\|_{L_T^\infty L_x^2(\mathbb{R}^n)},\\
&A_3 \lesssim \nu T  \sum_{\lambda=0, 1}\sum_{j=1}^n\sum_{k\in
\mathbb{Z}^n} \langle k\rangle^{2}\|\Box_k
 (\partial_{x_j}^\lambda
u)\|_{L_T^\infty L_x^2(\mathbb{R}^n)}.\label{4.94}
\end{align}
From the argument in \eqref{4.20}, and initial data belong to
$M_{2,1}^{4}$, we have $$\sum_{\lambda=0,
1}\sum_{j=1}^n\sum_{k\in \mathbb{Z}^n}\langle k\rangle^{3/2}\langle
k\rangle \|\Box_k
 (\partial_{x_j}^\lambda
u)\|_{L_t^\infty L_x^2(\mathbb{R}^\times \mathbb{R}^n)}\leq C.$$

Collection \eqref{4.93}--\eqref{4.94}, we finally obtain
\begin{align}
\|u_\nu-v\|_{Y_T}\lesssim  \nu T  \sum_{\lambda=0,
1}\sum_{j=1}^n\sum_{k\in \mathbb{Z}^n}\langle k\rangle^{3/2}\langle
k\rangle \|\Box_k
 (\partial_{x_j}^\lambda
u)\|_{L_t^\infty L_x^2(\mathbb{R}^\times\mathbb{R}^n)}.
\end{align}
In this way, we obtain the results of limit behavior
\begin{align}
\|u_\nu - v\|_{Y_T} \rightarrow 0, \quad \nu \rightarrow 0
.\nonumber
\end{align}

\section{Local well-posedness results for $n=1,2$}
When $n=1, 2$, $T\leq 1$,  define
\begin{align}
&\rho_1^T(u)=\sum_{i=1}^n  \sum_{k\in \mathbb{Z}^n, |k_i|>4}\langle
k_i\rangle^{s-1/2 } \|\Box_k u\|_{L_{x_i}^{\infty}L_{(x_j)j\neq i
}^{2}L_T^2( \mathbb{R}^n)}:=\sum_{i=1}^n \rho_1^i(u) \nonumber\\
&\rho_2^T(u)=\sum_{i=1}^n \sum_{k\in \mathbb{Z}^n} \|\Box_k
u\|_{L_{x_i}^{2}L_{(x_j)j\neq i }^{\infty}L_T^\infty(
\mathbb{R}^n)}\nonumber\\
&\rho_3^T(u)=\sum_{k\in \mathbb{Z}^n}\langle k\rangle^{s-1 }\|\Box_k
u\|_{L_T^\infty L_x^2( \mathbb{R}^n)}.\nonumber
\end{align}
 Define resolution space as following:
\begin{align}
X^T_s:= \{u\in \mathscr{S}'([0, T] \times \mathbb{R}^n):
\|u\|_{X^T_s}:= \sum_{l=1}^3\sum_{\lambda=0,
1}\sum_{j=1}^n\rho_l^T(\partial_{x_j}^\lambda u)\leq
\delta_0\}\label{xt}
\end{align}
We write ${X}^T_{5/2}$ as ${X}^T$ for short.

 Similar to the proof of global
well-posedness results,  we only need to consider the case $u_0 \in
M_{1, 1}^s, s=5/2$ is small enough and we have:
\begin{align}
\rho_2^T(\partial_{x_j}^\lambda G_\nu (t)u_0)\leq \sum_{i=1}^n
\sum_{k\in \mathbb{Z}^n}\langle k_i\rangle^{3/2} \|\Box_k
u_0\|_{L_1(\mathbb{R}^n)}\lesssim \|u_0\|_{M_{1, 1}^{3/2 }}\nonumber
\end{align}
$\rho_1^T(\partial_{x_j}^\lambda G_\nu (t)u_0),
\rho_3^T(\partial_{x_j}^\lambda G_\nu (t)u_0)$ are similar to
section 4. So we have
\begin{align}
\|G_\nu (t)u_0\|_{X^T} \lesssim \|u_0\|_{M_{1, 1}^{5/2 }}\label{5.1}
\end{align}
Notice for any $p\geq 1, q\geq 2$, we have
\begin{align}
\|\Box_k u\|_{L_T^{p}L_x^{q}( \mathbb{R}^n)} \lesssim
T^{1/p}\|\Box_k u\|_{L_T^\infty L_x^2( \mathbb{R}^n)}\label{6.2}
\end{align}
 From \eqref{3.24}, Lemma \ref{lem5.1} and H\"older's
inequality we have
\begin{align}
&\rho_2^T(\partial_{x_j}^\lambda(\mathscr{A}(\lambda_1^i(\partial_{x_i}\bar{u})
u^2))) \nonumber\\
&\lesssim  \sum_{l=1}^n\sum_{k\in \mathbb{Z}^n} \langle
k_j\rangle^{\lambda}\langle k_l\rangle^{1/2}\|\Box_k (\partial_{x_i}\bar{u})
u^2\|_{L_{T}^1 L_x^1(
\mathbb{R}^n)}\nonumber\\
&\lesssim n \Big[ \sum_{k\in \mathbb{Z}^n} \langle
k\rangle^{3/2 }\|\Box_k
(\partial_{x_i}\bar{u})\|_{L_{T}^3 L_x^3( \mathbb{R}^n)}\times
\big(\sum_{k\in \mathbb{Z}^n} \langle
k\rangle^{3/2 }\|\Box_k u\|_{L_{T}^3 L_x^3(
\mathbb{R}^n)}\big)^2\Big]\nonumber\\
&\lesssim n T \rho_3(\partial_{x_i}u)\rho_3(u)^2\label{5.4}
\end{align}
From Proposition \ref{4.12} and Proposition \ref{p2.9.1} and the
estimate in \eqref{4.74}, \eqref{5.3}--\eqref{5.6}, we have
\begin{align}
&\rho_3^T(\partial_{x_j}^\lambda(\mathscr{A}(\lambda_1^i(\partial_{x_i}\bar{u})
u^2)))\nonumber\\
&\lesssim n\sum_{k\in \mathbb{Z}^n, |k_j|>4}\langle
 k_j\rangle^{\lambda+3/2-1/2}
\sum_{k^{(1)}, k^{(2)}, k^{(3)}}\|\Box_k (\Box_{k^{(1)}}
\partial_{x_i}\bar{u}
\Box_{k^{(2)}}u\Box_{k^{(3)}}u)\|_{L_{x_j}^{1}L_{x_i(i\neq
j)}^{2}L_T^2(
\mathbb{R}^n)}\nonumber\\
&+ n\sum_{k\in \mathbb{Z}^n, |k_j|\leq 4}\langle
 k_j\rangle^{5/2}
\sum_{k^{(1)}, k^{(2)}, k^{(3)}}\|\Box_k (\Box_{k^{(1)}}
\partial_{x_i}\bar{u}
\Box_{k^{(2)}}u\Box_{k^{(3)}}u)\|_{L^{4/3}_T
L^{4/3}_x(\mathbb{R}^n)}\nonumber\\
&\lesssim \rho_1^T(\partial_{x_i}u)\rho_2^T(u)^2 +
T\rho_3^T(\partial_{x_i}u)\rho_3^T(u)^2.\label{6.5}
\end{align}
$\rho_1^T(\partial_{x_j}^\lambda(\mathscr{A}(\lambda_1^i(\partial_{x_i}\bar{u})
u^2)))$ is similar to section 4, we omit the detail. The estimate
for $\|\mathscr{A}(\alpha |u|^{2\delta}u)\|_{X^T_s}$ is similar, we
do not repeat  here.

From the above, we obtain for any $T\leq 1,$
\begin{align}
\|u\|_{X^T}\leq (1+T)\|u_0\|_{M_{1, 1}^{5/2
 }} + (1+T)(\|u\|_{X^T}^3 + \|u\|_{X^T}^{2\delta +1}).
\end{align}
Using the small initial data which is independent of $T$, we have
$u\in X_T$ and $\|u\|_{X_T}$ small enough. Similar to \eqref{4.20},
we obtain that
\begin{align}
\|u\|_{X^T_s}\leq \|u_0\|_{M_{1, 1}^{s
 }} + (1+T)\|u\|_{X^T_s}(\|u\|_{X^T}^2 + \|u\|_{X^T}^{2\delta}),
 \quad s> 5/2\nonumber.
\end{align}

Now, we show that $u\in L^\infty_T (M_{1, 1}^{s-1/2}
(\mathbb{R}^n))$. Similar to the estimate of $\rho_3$ in Section 4,
we also have
\begin{align}
&\sum_{k\in \mathbb{Z}^n}\langle
 k\rangle^{s-1/2}\|\Box_k [\overrightarrow{\lambda_1}\cdot\nabla(|u|^2u) +
(\overrightarrow{\lambda_2}\cdot\nabla u)|u|^2]\|_{L^1_TL_x^1(
 \mathbb{R}^n)}\nonumber\\
& \lesssim  \sum_{j=1}^n
 \sum_{k\in \mathbb{Z}^n, |k_j|>4}\langle
 k_j\rangle^{s-1/2}
\sum_{\mathbb{S}_{2, 1}^{(j)}}\|\Box_k (\Box_{k^{(1)}}
\partial_{x_i}\bar{u}
\Box_{k^{(2)}}u\Box_{k^{(3)}}u)\|_{L^1_TL_x^1(
 \mathbb{R}^n)}\nonumber\\
&+ \sum_{j=1}^n
 \sum_{k\in \mathbb{Z}^n, |k_j|\geq  4}\langle
 k_j\rangle^{s-1/2}
\sum_{\mathbb{S}_{2, 2}^{(j)}}\|\Box_k (\Box_{k^{(1)}}
\partial_{x_i}\bar{u}
\Box_{k^{(2)}}u\Box_{k^{(3)}}u)\|_{L^1_TL_x^1(
 \mathbb{R}^n)}\nonumber\\
 &+ \sum_{j=1}^n
 \sum_{k\in \mathbb{Z}^n, |k_j|\leq 4}\langle
 k_j\rangle^{s-1/2}
\sum_{k^{(1)}, k^{(2)}, k^{(3)}}\|\Box_{k^{(1)}}
\partial_{x_i}\bar{u}
\Box_{k^{(2)}}u\Box_{k^{(3)}}u\|_{L^1_TL_x^1(
 \mathbb{R}^n)}.\nonumber\\
 &\leq \rho_1^T({\partial_{x_i}u})(\rho_2^T(u) +
\rho_3^T(u))^2+ \rho_3^T(\partial_{x_i}u)\rho_3^T(u)^2, \label{6.7}
\end{align}
where we use
\begin{align}
&\|\Box_{k^{(1)}} (\partial_{x_i}\bar{u})
\Box_{k^{(2)}}u\Box_{k^{(3)}}u\|_{L_T^{1}L_x^1(
\mathbb{R}^n)}\nonumber\\
&\leq \|\Box_{k^{(1)}} (\partial_{x_i}\bar{u})
|\Box_{k^{(2)}}u\Box_{k^{(3)}}u|^{1/2}\|_{L_T^{2}L_x^2(
\mathbb{R}^n)}\|
|\Box_{k^{(2)}}u\Box_{k^{(3)}}u|^{1/2}\|_{L_T^{2}L_x^2(
\mathbb{R}^n)}\nonumber\\
&\lesssim T^{1/2}\|\Box_{k^{(1)}} (\partial_{x_i}\bar{u})\|
_{L_{x_j}^{\infty}L_{x_r(r\neq j)}^{2}L_T^2( \mathbb{R}^n)}
\|\Box_{k^{(2)}}u \|^{1/2}_{L_{x_j}^{1}L_{x_r(r\neq
j)}^{\infty}L_T^\infty( \mathbb{R}^n)} \| \Box_{k^{(3)}}u
\|^{1/2}_{L_{x_j}^{1}L_{x_r(r\neq
j)}^{\infty}L_T^\infty(\mathbb{R}^n)}\nonumber\\
&\quad\quad\times \|\Box_{k^{(2)}} u\|^{1/2}_{L_T^\infty
L_x^2(\mathbb{R}^n)} \| \Box_{k^{(3)}}u \|^{1/2}_{L_T^\infty
L_x^2(\mathbb{R}^n)}\nonumber
\end{align}
and
\begin{align}
\|\Box_k u\|_{L_{x_1}^{2}L_{x_2, \ldots, x_n}^{\infty}L_T^\infty(
\mathbb{R}^n)}\leq \|\Box_k u\|_{L_{x_1}^{1}L_{x_2, \ldots,
x_n}^{\infty}L_T^\infty( \mathbb{R}^n)}.
\end{align}
 From Proposition \ref{ubst} and Minkowski's inequality, we have
\begin{align}
&\sum_{k\in \mathbb{Z}^n}\langle
 k\rangle^{s-1/2}\|\Box_k u\|_{L^\infty_TL_x^1(
 \mathbb{R}^n)}\nonumber\\
&\leq \sum_{k\in \mathbb{Z}^n}\langle
 k\rangle^{s-1/2}\|\Box_k G_\nu(t) u_0\|_{L^\infty_TL_x^1(\mathbb{R}^+ \times
 \mathbb{R}^n)} \nonumber\\
 & +\sum_{k\in \mathbb{Z}^n}\langle
 k\rangle^{s-1/2}\|\Box_k \mathscr{A}[\overrightarrow{\lambda_1}\cdot\nabla(|u|^2u) +
(\overrightarrow{\lambda_2}\cdot\nabla u)|u|^2+ \alpha
|u|^{2\delta}u]\|_{L^\infty_TL_x^1(
 \mathbb{R}^n)}\nonumber\\
&\leq \|u_0\|_{M_{2, 1}^{s-1/2}} + \sum_{k\in \mathbb{Z}^n}\langle
 k\rangle^{s-1/2}\|\Box_k [\overrightarrow{\lambda_1}\cdot\nabla(|u|^2u) +
(\overrightarrow{\lambda_2}\cdot\nabla u)|u|^2+ \alpha
|u|^{2\delta}u]\|_{L^1_TL_x^1(
 \mathbb{R}^n)}\nonumber\\
&\lesssim \|u_0\|_{M_{2, 1}^{s-1/2}} +  \|u\|_{X_s}^3+
\|u\|_{X_s}^{2\delta+1}.\label{6.8}
\end{align}
We obtain the local well-posedness results. The limit behavior
results are almost the same as in Section 5.

\section{Local well-posedness for the quadratic DNLS}
In this section, we will prove local well-posedness results for
equation
\begin{align}
u_t=  (\nu + {\rm i})\triangle u + \vec{\lambda}\cdot \nabla(u^2), \quad
u(0,x)=u_0(x),\label{GL2}
\end{align}
and equation
\begin{align}
v_t=  {\rm }i\triangle v + \vec{\lambda}\cdot \nabla(v^2), \quad
v(0,x)=v_0(x). \label{s2}
\end{align}
When $n \in \mathbb{N}$, $T\leq 1$,  define
\begin{align}
&\rho_1^T(u)=\sum_{i=1}^n  \sum_{k\in \mathbb{Z}^n, |k_i|>4}\langle
k_i\rangle^{s- 1/2 } \|\Box_k u\|_{L_{x_i}^{\infty}L_{(x_j)j\neq i
}^{2}L_T^2( \mathbb{R}^n)}:=\sum_{i=1}^n \rho_1^{i}(u), \nonumber\\
&\rho_2^T(u)=\sum_{i=1}^n \sum_{k\in \mathbb{Z}^n} \|\Box_k
u\|_{L_{x_i}^{1}L_{(x_j)j\neq i }^{\infty}L_T^\infty(
\mathbb{R}^n)},\nonumber\\
&\rho_3^T(u)=\sum_{k\in \mathbb{Z}^n}\langle k\rangle^{s- 1}\|\Box_k
u\|_{L_T^\infty L_x^2 ( \mathbb{R}^n)},\nonumber
\end{align}
Define resolution space as following:
\begin{align}
{{\tilde{X}}^T_s:= \{u\in \mathscr{S}'([0, T] \times \mathbb{R}^n):
\|u\|_{\tilde{X}^T_s}}:= \sum_{l=1}^3\sum_{j=1}^n\rho_l^T( u)\leq
\delta_0\}.\label{xt1}
\end{align}
We write $\tilde{X}^T_{3}$ as $\tilde{X}^T$ for short.

 We solve
equation \eqref{GL2} first.
 Similar to the proof of global
well-posedness results,  we only need to consider the case $u_0 \in
M_{1, 1}^s, s=3$ is small enough.
  Similar to section 4, We have
\begin{align}
\|G_\nu (t)u_0\|_{\tilde{X}^T} \lesssim \|u_0\|_{M_{1, 1}^{3
 }}.\label{5.1}
\end{align}
From \eqref{3.24}, Lemma \ref{lem5.1} and H\"older's inequality we
have
\begin{align}
\rho_2^T(\mathscr{A}(\partial_{x_j}u^2) ) &\lesssim
\sum_{l=1}^n\sum_{k\in \mathbb{Z}^n} \langle k_j\rangle \langle
k_l\rangle^{
 }\|\Box_k u^2\|_{L_{T}^1 L_x^1(
\mathbb{R}^n)}\nonumber\\
&\lesssim n \Big[ \sum_{k\in \mathbb{Z}^n} \langle k\rangle^{2
 }\|\Box_k u\|_{L_{T}^2 L_x^2( \mathbb{R}^n)}\Big]^2\nonumber\\
&\lesssim n T \rho_3^T(u)^2.\label{}
\end{align}
Similar to \eqref{5.3},
\begin{align}
&\rho_1^{1}(\mathscr{A}(\partial_{x_1}u^2) )\nonumber\\
 &\lesssim \sum_{k\in \mathbb{Z}^n, |k_1|>4}\langle
 k_1\rangle^{5/2}
\sum_{\mathbb{S}_{2, 1}^{(1)}}\|\Box_k (\Box_{k^{(1)}} u
\Box_{k^{(2)}}u)\|_{L_{x_1}^{1}L_{x_2, \ldots,
x_n}^{2}L_T^2( \mathbb{R}^n)}\nonumber\\
&\quad +\sum_{k\in \mathbb{Z}^n, |k_1|\leq 4}\langle
 k_1\rangle^{7/2 }
\sum_{\mathbb{S}_{2, 2}^{(1)}}\|\Box_k (
\Box_{k^{(1)}}u\Box_{k^{(2)}}u)\|_{L^{1}_T L^{2}_x
(\mathbb{R}^n)}\nonumber\\
&:= I + II\label{}.
\end{align}
 Notice \eqref{6.2} and
\begin{align}
\|\Box_k u\|_{L_{x_1}^{2}L_{x_2, \ldots, x_n}^{\infty}L_T^\infty(
\mathbb{R}^n)}\leq \|\Box_k u\|_{L_{x_1}^{1}L_{x_2, \ldots,
x_n}^{\infty}L_T^\infty( \mathbb{R}^n)}.
\end{align}
Similar to \eqref{5.5}, \eqref{5.6}, we have
\begin{align}
I \leq \rho_1^{T}(u) \rho_2^{T}(u),\quad \quad II\leq
\rho_2^{T}(u)^2.
\end{align}
So, $\rho_1^{T}(\mathscr{A}(\partial_{x_j}u^2) )\leq \rho_1^{T}(u)
\rho_2^{T}(u)+ \rho_2^{T}(u)^2$. We estimate
$\rho_3^T(\mathscr{A}(\partial_{x_j}u^2) )$ via a similar way as in
\eqref{6.5}:
\begin{align}
&\rho_3^T(\mathscr{A}(\partial_{x_j}u^2) )\nonumber\\
&\lesssim n\sum_{k\in \mathbb{Z}^n, |k_j|>4}\langle
 k_j\rangle^{3}
\sum_{k^{(1)}, k^{(2)}, k^{(3)}}\|\Box_k (\Box_{k^{(1)}} u
\Box_{k^{(2)}}u)\|_{L_{x_j}^{1}L_{x_i(i\neq j)}^{2}L_T^2(
\mathbb{R}^n)}\nonumber\\
&+ n\sum_{k\in \mathbb{Z}^n, |k_j|\leq 4}\langle
 k_j\rangle^{7/2}
\sum_{k^{(1)}, k^{(2)}, k^{(3)}}\|\Box_k (\Box_{k^{(1)}} u
\Box_{k^{(2)}}u)\|_{L^{4/3}_T
L^{4/3}_x(\mathbb{R}^n)}\nonumber\\
&\lesssim \rho_1^T(u)\rho_2^T(u) + T\rho_3^T(u)^2.
\end{align}
So we obtain
\begin{align}
\|u\|_{\tilde{X}^T_3}\leq \|u_0\|_{M_{1, 1}^3}+(1+T)
\|u\|^2_{\tilde{X}^T_3},\nonumber
\end{align}
$\|u\|_{\tilde{X}^T_3}$ is sufficiently small, and also
\begin{align}
\|u\|_{\tilde{X}^T_s}\leq \|u_0\|_{M_{1, 1}^s}+ (1+T)
\|u\|_{\tilde{X}^T_s}\|u\|_{\tilde{X}^T_3}.\nonumber
\end{align}
 In this way, we can also obtain local well-posedness
results for the solution $v$ of Schr\"odinger equation \eqref{s2}.
The inviscid limit for \eqref{GL2} is almost the same as section 5.
We can obtain
\begin{align}
\|u_\nu-v\|_{\tilde{X}^T_3}\rightarrow 0, \quad \nu \rightarrow 0.
\end{align}
We omit the detail here.

Via a similar way as in \eqref{6.7} and \eqref{6.8}, we show that
$u\in L^\infty_T (M_{1, 1}^{s-1/2} (\mathbb{R}^n))$, $s>3$. where we
use
\begin{align}
&\|\Box_k (\Box_{k^{(1)}} u \Box_{k^{(2)}}u)\|_{L_T^{1}L_x^1(
\mathbb{R}^n)}\nonumber\\
&\leq \|\Box_{k^{(1)}}u |\Box_{k^{(2)}}u|^{1/2}\|_{L_T^{2}L_x^2(
\mathbb{R}^n)}\| |\Box_{k^{(2)}}u|^{1/2}\|_{L_T^{2}L_x^2(
\mathbb{R}^n)}\nonumber\\
&\lesssim T^{1/2}\|\Box_{k^{(1)}} u\| _{L_{x_j}^{\infty}L_{x_r(r\neq
j)}^{2}L_T^2( \mathbb{R}^n)} \|\Box_{k^{(2)}}u
\|^{1/2}_{L_{x_j}^{1}L_{x_r(r\neq j)}^{\infty}L_T^\infty(
\mathbb{R}^n)} \|\Box_{k^{(2)}} u\|^{1/2}_{L_T^\infty
L_x^1(\mathbb{R}^n)} .\nonumber
\end{align}
Then as  the estimate of \eqref{6.7} and \eqref{6.8} in  Section 6,
we have
\begin{align}
\|u\|_{K} \lesssim \|u_0\|_{M_{1, 1}^{s-1/2}} +
\|u\|_{_{\tilde{X}^T_s}}\|u\|_{_{\tilde{X}^T_3}}^{1/2}\|u\|_{K}^{1/2},
\end{align}
where
\begin{align}
\|u\|_K:= &\sum_{k\in \mathbb{Z}^n}\langle
 k\rangle^{s-1/2}\|\Box_k u\|_{L^\infty_TL_x^1(\mathbb{R}^+ \times
 \mathbb{R}^n)}.\nonumber\\
\end{align}
Using $\|u\|_{\tilde{X}^T_3}$ is sufficiently small, we obtain $u\in
K$ and so $u\in L^\infty_T (M_{1, 1}^{s-1/2} (\mathbb{R}^n))$.

\newpage

\begin{appendix}

\section{Appendix}

\subsection{A nonlinear estimate}
In this section, we will show the detail proof of \eqref{5.11}. We
follow the idea in our earlier work \cite{WH}. From \eqref{1.22},
\eqref{2.32}, Lemma \ref{lem2.4} and \eqref{5.2}, (where we put
$\ell=2$), obtained:

\begin{align}
&\rho_1^1(\partial_{x_j}^{\lambda}(\mathscr{A}\lambda_1^i(\partial_{x_i}\bar{u})
u^2))\nonumber\\
 &\lesssim \sum_{k\in \mathbb{Z}^n, |k_1|>4}\langle
 k_1\rangle^{\lambda+3/2}
\sum_{\mathbb{S}_{2, 1}^{(1)}}\|\Box_k (\Box_{k^{(1)}}
\partial_{x_i}\bar{u}
\Box_{k^{(2)}}u\Box_{k^{(3)}}u)\|_{L_{x_1}^{1}L_{x_2, \ldots,
x_n}^{2}L_t^2(\mathbb{R}_+
\times \mathbb{R}^n)}\nonumber\\
&+\sum_{k\in \mathbb{Z}^n, |k_1|>4}\langle
 k_1\rangle^{\lambda+2}
\sum_{\mathbb{S}_{2, 2}^{(1)}}\|\Box_k (\Box_{k^{(1)}}
\partial_{x_i}\bar{u}
\Box_{k^{(2)}}u\Box_{k^{(3)}}u)\|_{L^{1}_t(\mathbb{R}^+;
L^{2}_x(\mathbb{R}^n))}\nonumber\\
&:= I + II\label{5.3}
\end{align}
In view of the support property of $\widehat{\Box_k u}$, we can see
\begin{align}
\Box_k (\Box_{k^{(1)}} (\partial_{x_i}\bar{u})
\Box_{k^{(2)}}u\Box_{k^{(3)}}u)=0, \quad if  \quad
|k-k^{(1)}-k^{(2)}-k^{(3)}|\geq C.\nonumber
\end{align}
 For $I$, since $|k-k^{(1)}-k^{(2)}-k^{(3)}| \leq C$, it is easy to see $|k_1|\leq C
 \max_{r=1,2,3}|k_1^r|$, we can assume $
 |k_1^{(1)}|=\max_{r=1,2,3}|k_1^r|$.
 From H\"older's inequality
\begin{align}
&I = \sum_{k\in \mathbb{Z}^n, |k_1|>4}\langle
 k_1\rangle^{5/2}
\sum_{\mathbb{S}_{2, 1}^{(1)}}\|(\Box_{k^{(1)}}
(\partial_{x_i}\bar{u})
\Box_{k^{(2)}}u\Box_{k^{(3)}}u)\|_{L_{x_1}^{1}L_{x_2, \ldots,
x_n}^{2}L_t^2(\mathbb{R}_+
\times \mathbb{R}^n)}\nonumber\\
&\lesssim \sum_{k^{(1)}\in \mathbb{Z}^n, |k_1^{(1)}|>4}\langle
k_1^{(1)}\rangle^{5/2} \|\Box_{k^{(1)}}
\partial_{x_i}\bar{u}\|_{L_{x_1}^{\infty}L_{x_2, \ldots,
x_n}^{2}L_t^2(\mathbb{R}_+
\times \mathbb{R}^n)}\nonumber\\
&\times  \Big(\sum_{k^{(2)}\in \mathbb{Z}^n}\|\Box_{k^{(2)}}
u\|_{L_{x_1}^{2}L_{x_2, \ldots, x_n}^{\infty}L_t^\infty(\mathbb{R}_+
\times \mathbb{R}^n)}\Big)^2\nonumber\\
&\lesssim \rho_1({\partial_{x_i}u})\rho_2(u)^2.\label{5.5}
\end{align}

For $II$, from the definition of $\mathbb{S}_{\ell, 2}^{(1)}$, and
 $|k-k^{(1)}-k^{(2)}-k^{(3)}| \leq C$, we have $|k_1|\leq C$. If we fix $k^{(1)}, k^{(2)}, k^{(3)}$,
 then $k$ in the sum is finite.
 So from Holder inequality we have:
\begin{align}
 II &\lesssim  \sum_{\mathbb{S}_{2, 2}^{(1)}}\|\Box_{k^{(1)}}
\partial_{x_i}\bar{u}
\Box_{k^{(2)}}u\Box_{k^{(3)}}u\|_{L^{1}_t(\mathbb{R}^+;
L^{2}_x(\mathbb{R}^n))}\nonumber\\
&\lesssim \sum_{k^{(1)}\in \mathbb{Z}^n}\|\Box_{k^{(1)}}
\partial_{x_i}\bar{u}
\|_{ L_t^3L_x^6(\mathbb{R}^+ \times
\mathbb{R}^n)}\Big(\sum_{k^{(2)}\in \mathbb{Z}^n}\|\Box_{k^{(2)}} u
\|_{ L_t^3L_x^6(\mathbb{R}^+ \times
\mathbb{R}^n)}\Big)^2\nonumber\\
&\lesssim \rho_3(\partial_{x_i}u) \rho_3(u)^2.\label{5.6}
\end{align}

Now we estimate
$\rho_1^1(\partial_{x_2}^\lambda\mathscr{A}(\lambda_1^i(\partial_{x_i}\bar{u})
u^2)), \lambda=1$ (when $\lambda=0$ then it is the same to the case
$j=1$).
 Let $P_i:= \mathscr{F}_{\xi_1,
\xi_2}^{-1}\psi_i \mathscr{F}_{x_1, x_2}$, where $\psi_i (i=1, 2)$
be as in Lemma \ref{lem4.2}. From  Lemma \ref{lem4.2}, we have
\begin{align}
&\rho_1^1(\partial_{x_2}\mathscr{A}(\lambda_1^i(\partial_{x_i}\bar{u})
u^2) \nonumber\\
&\lesssim \sum_{k\in \mathbb{Z}^n, |k_1|>4}\langle
 k_1\rangle^2
\|P_1 \partial_{x_2}\Box_k
(\mathscr{A}\lambda_1^i(\partial_{x_i}\bar{u})
u^2)\|_{L_{x_1}^{\infty}L_{x_2, \ldots, x_n}^{2}L_t^2(\mathbb{R}_+
\times \mathbb{R}^n)}\nonumber\\
&+ \sum_{k\in \mathbb{Z}^n, |k_1|>4}\langle
 k_1\rangle^2
\|P_2 \partial_{x_2}\Box_k
(\mathscr{A}\lambda_1^i(\partial_{x_i}\bar{u})
u^2)\|_{L_{x_1}^{\infty}L_{x_2, \ldots, x_n}^{2}L_t^2(\mathbb{R}_+
\times \mathbb{R}^n)}\nonumber\\
&:= III + IV.
\end{align}

For $III$, applying decomposition \eqref{5.2}(where we consider
variable $x_1$, namely, divide into $\mathbb{S}_{2, 1}^{(1)}$ and
~$\mathbb{S}_{2, 2}^{(1)}$ in \eqref{5.2}), then from \eqref{4.5},
and \eqref{3.2}, we have
\begin{align}
III\lesssim  & \sum_{k\in \mathbb{Z}^n, |k_1|>4}\langle
 k_1\rangle^{5/2}
\|P_1 \partial_{x_2}\Box_k \mathscr{A}\sum_{\mathbb{S}_{2,
1}^{(1)}}(\Box_{k^{(1)}}
\partial_{x_i}\bar{u}
\Box_{k^{(2)}}u\Box_{k^{(3)}}u)\|_{L_{x_1}^{\infty}L_{x_2, \ldots,
x_n}^{2}L_t^2(\mathbb{R}^+
\times \mathbb{R}^n)}\nonumber\\
&+\sum_{k\in \mathbb{Z}^n, |k_1|>4,  |k_2|\lesssim |k_1|}\langle
 k_1\rangle^{5/2}
\|P_1 \partial_{x_2}\Box_k \mathscr{A}\sum_{\mathbb{S}_{2,
2}^{(1)}}(\Box_{k^{(1)}}
\partial_{x_i}\bar{u}
\Box_{k^{(2)}}u\Box_{k^{(3)}}u)\|_{L_{x_1}^{\infty}L_{x_2, \ldots,
x_n}^{2}L_t^2(\mathbb{R}^+
\times \mathbb{R}^n)}\nonumber\\
& \lesssim \sum_{k\in \mathbb{Z}^n, |k_1|>4}\langle
 k_1\rangle^{5/2}
\sum_{\mathbb{S}_{2, 1}^{(1)}}\|\Box_k (\Box_{k^{(1)}}
\partial_{x_i}\bar{u}
\Box_{k^{(2)}}u\Box_{k^{(3)}}u)\|_{L_{x_1}^{1}L_{x_2, \ldots,
x_n}^{2}L_t^2(\mathbb{R}^+
\times \mathbb{R}^n)}\nonumber\\
&+\sum_{k\in \mathbb{Z}^n, |k_1|>4, |k_2|\lesssim |k_1|}\langle
 k_1\rangle^{2}\langle k_2\rangle
\sum_{\mathbb{S}_{2, 2}^{(1)}}\|\Box_k (\Box_{k^{(1)}}
\partial_{x_i}\bar{u}
\Box_{k^{(2)}}u\Box_{k^{(3)}}u)\|_{L^{1}_t(\mathbb{R}^+;
L^{2}_x(\mathbb{R}^n))}\nonumber\\
&\lesssim
\rho_1^1({\partial_{x_i}u})\rho_2(u)^2+\rho_3(\partial_{x_i}u)
\rho_3(u)^2.\label{5.9}
\end{align}
At the last step of \eqref{5.9}, notice the definition of
$\mathbb{S}_{2, 2}^{(1)}$, it is easy to see $ |k_2|\lesssim
|k_1|\leq C$, then it comes back to \eqref{5.5}, \eqref{5.6}, so
repeat the proof of \eqref{5.5}, \eqref{5.6}, we can obtain
\eqref{5.9}, as desired.

For the estimates of $IV$, applying the decomposition \eqref{5.2}
(where we consider variable $x_2$, namely, divide into
$\mathbb{S}_{2, 1}^{(2)}$ and $\mathbb{S}_{2, 2}^{(2)}$ in
\eqref{5.2}), in addition to \eqref{4.6} and \eqref{3.2}, we have
\begin{align}
IV\lesssim  & \sum_{k\in \mathbb{Z}^n, |k_1|>4}\langle
 k_1\rangle^{5/2}
\|P_2 \partial_{x_2}\Box_k \mathscr{A}\sum_{\mathbb{S}_{2,
1}^{(2)}}(\Box_{k^{(1)}}
\partial_{x_i}\bar{u}
\Box_{k^{(2)}}u\Box_{k^{(3)}}u)\|_{L_{x_1}^{\infty}L_{x_2, \ldots,
x_n}^{2}L_t^2(\mathbb{R}^+
\times \mathbb{R}^n)}\nonumber\\
&+\sum_{k\in \mathbb{Z}^n, |k_1|>4,  |k_1|\lesssim |k_2|}\langle
 k_1\rangle^{5/2}
\|P_2 \partial_{x_2}\Box_k \mathscr{A}\sum_{\mathbb{S}_{2,
2}^{(2)}}(\Box_{k^{(1)}}
\partial_{x_i}\bar{u}
\Box_{k^{(2)}}u\Box_{k^{(3)}}u)\|_{L_{x_1}^{\infty}L_{x_2, \ldots,
x_n}^{2}L_t^2(\mathbb{R}^+
\times \mathbb{R}^n)}\nonumber\\
& \lesssim \sum_{k\in \mathbb{Z}^n, |k_2|>4}\langle
 k_2\rangle^{5/2}
\sum_{\mathbb{S}_{2, 1}^{(2)}}\|\Box_k (\Box_{k^{(1)}}
\partial_{x_i}\bar{u}
\Box_{k^{(2)}}u\Box_{k^{(3)}}u)\|_{L_{x_1}^{1}L_{x_2, \ldots,
x_n}^{2}L_t^2(\mathbb{R}^+
\times \mathbb{R}^n)}\nonumber\\
&+\sum_{k\in \mathbb{Z}^n, |k_1|>4, |k_1|\lesssim |k_2|}\langle
 k_1\rangle^{2}\langle
 k_2\rangle
\sum_{\mathbb{S}_{2, 2}^{(2)}}\|\Box_k (\Box_{k^{(1)}}
\partial_{x_i}\bar{u}
\Box_{k^{(2)}}u\Box_{k^{(3)}}u)\|_{L^{1}_t(\mathbb{R}^+;
L^{2}_x(\mathbb{R}^n))}\nonumber\\
 &\lesssim
\rho_1^2({\partial_{x_i}u})\rho_2(u)^2+\rho_3(\partial_{x_i}u)
\rho_3(u)^2.\label{5.10}
\end{align}
At the last step of \eqref{5.10}, notice the definition of
$\mathbb{S}_{2, 2}^{(2)}$, it is easy to see $ |k_1|\lesssim
|k_2|\leq C$. This way, it comes back to \eqref{5.5}, \eqref{5.6},
follow the same process there, we can obtain \eqref{5.10}, as
desired.

\subsection{Appendix}

\begin{lem}\label{appb}
For any $s\in \mathbb{R}$ and any $s^+ >s$, there exist $ \theta >0$
such that
\begin{align}
\|f\|_{M_{2, 1}^s}\leq C^\theta \|f\|_{M_{2,
1}^{s^+}}^{1-2\theta}\|f\|_{L^2}^\theta,
\end{align}
where $s^+ = \frac{s+ 2\theta}{1-2\theta}$.
\end{lem}
\noindent{\bf Proof:}
\begin{align}
\|f\|_{M_{2, 1}^s}&=\sum_{k \in \mathbb{Z}^n} \|\Box_k
f\|_{H^s}\nonumber\\
&\leq \sum_{k \in \mathbb{Z}^n} \|\Box_k
f\|_{H^{s+\epsilon}}^{1-\theta}\|\Box_k f\|_{L^2}^\theta\nonumber\\
&\leq \sum_{k \in \mathbb{Z}^n}\langle
k\rangle^{(s+\epsilon)(1-\theta)}\|\Box_k
f\|_{L^2}^{1-2\theta}\|\Box_k f\|_{L^2}^{2\theta}\nonumber\\
&\leq \Big\|\langle k\rangle^{(s+\epsilon)(1-\theta)}\|\Box_k
f\|_{L^2}^{1-2\theta}\Big\|_{l^{\frac{1}{1-\theta}}}\Big\|\|\Box_k
f\|_{L^2}^{2\theta}\Big\|_{l^{\frac{1}{\theta}}}\nonumber\\
&\leq \Big\|\langle
k\rangle^{(s+\epsilon+\frac{2\theta}{1-\theta})}\|\Box_k
f\|_{L^2}^{\frac{1-2\theta}{1-\theta}}\Big\|_{l^{\frac{1-\theta}{1-2\theta}}}^{1-\theta}
\Big\|\langle k
\rangle^{-\frac{2\theta}{1-\theta}}\Big\|_{l^{\frac{1-\theta}{\theta}}}^{1-\theta}
\Big\|f\Big\|_{L^2}^\theta\nonumber\\
&\leq C^\theta \|f\|_{M_{2,
1}^{s^+}}^{1-2\theta}\|f\|_{L^2}^\theta\nonumber
\end{align}
where $\theta = \frac{\epsilon}{s+\epsilon}$, $s^+ =
\frac{1-\theta}{1-2\theta}(s+\epsilon+\frac{2\theta}{1-\theta})$.
$\hfill\Box$

\begin{lem}\label{appb1}
For any $s\in \mathbb{R}$ and any $s^+ >s$, there exist $ \theta >0$
such that
\begin{align}
\|f\|_{M_{1, 1}^s}\leq C  \|f\|_{M_{1,
1}^{s^+}}^{1-\theta}\|f\|_{L^1}^\theta,
\end{align}
where $s^+ = \frac{s+2\theta}{1-\theta}$.
\end{lem}
\noindent{\bf Proof:}
\begin{align}
\|f\|_{M_{2, 1}^s}&=\sum_{k \in \mathbb{Z}^n} \langle
k\rangle^{s}\|\Box_k
f\|_{L^1}\nonumber\\
&\leq \sum_{k \in \mathbb{Z}^n}\langle k\rangle^{s}\|\Box_k
f\|_{L^1}^{1-\theta}\|\Box_k f\|_{L^1}^{\theta}\nonumber\\
&\leq \sup_{k}\|\Box_k f\|_{L^1}^{\theta}\sum_{k \in
\mathbb{Z}^n}\langle k\rangle^{s}\|\Box_k
f\|_{L^1}^{1-\theta}\nonumber\\
 &\leq \|f\|_{L^1}^{\theta}\Big\|\langle
k\rangle^{(s+2\theta)}\|\Box_k
f\|_{L^1}^{1-\theta}\Big\|_{l^{\frac{1}{1-\theta}}}\Big\|\langle k
\rangle^{-2\theta}\Big\|_{l^{\frac{1}{\theta}}}\nonumber\\
&\leq C^\theta\|f\|_{L^1}^{\theta}\Big\|\langle
k\rangle^{\frac{s+2\theta}{1-\theta}}\|\Box_k
f\|_{L^1}\Big\|_{l^{1}}^{1-\theta}\nonumber\\
&\leq C^\theta \|f\|_{M_{1,
1}^{s^+}}^{1-\theta}\|f\|_{L^1}^\theta,\nonumber
\end{align}
where $s^+ = \frac{s+2\theta}{1-\theta}$.
\end{appendix}

\noindent{\bf Acknwoledgments}

The author is grateful to Z. Guo for his enlightening advice.


\medskip
\footnotesize
\begin{thebibliography}{31}

\bibitem{BJL} J. Bergh, L\"ofstr\"om, J. Interpolation spaces.
Grundlehren der Mathematischen Wissenschaften, 223. Springer,
Berlin-New York, 1976.

\bibitem{BJ} P. Bechouche, A. J\"{u}ngel,  Inviscid limits of the complex
Ginzburg--Landau equation, Commun. Math. Phys.,  {\bf 214} (2000),
201--226.

\bibitem{2} H. R. Brand, R. J. Deissler, Interaction of localized solutions
for subcritical bifurcations, Phys. Rev. Lett.,  {\bf 63} (1989)
2801--2804.

\bibitem{Ch} Micheal Christ, Ill posedness of a Schr\"odinger equation with
derivative nonlinearity. Preprint.

\bibitem{CT}  T. Cazenave,  F. B. Weissler, The cauchy problem for the
critical nonlinear Schr\"odingger equation in $H^s$. Nonlinear
Anal., {\bf 14 } (1990),   807--836.

\bibitem{4} M. C. Cross, P. C. Hohenberg, Pattern formation outside of
equilibrium, Rev. Mod. Phys., {\bf 65} (1993) 851--1089.

\bibitem{Chiha}  H. Chihara, Gain of regularity for
semilinear Schr\"odinger equations, Math. Ann., {\bf 315} (1999),
529--567.


\bibitem{JDE}  J. Duan, P. Holmes, and E. S. Titi, Regularity
approximation and asymptotic dynamics for a generalized
Ginzburg--Landau equation, Nonlinearity., {\bf 6} (1993), 915--933.


\bibitem{Fei2} H. G. Feichtinger, Modulation spaces on locally
compact Abelian group, Technical Report, University of Vienna, 1983.
Published in: ``Proc. Internat. Conf. on Wavelet and Applications",
99--140.  New Delhi Allied Publishers, India, 2003.
{http://www.unive.ac.at/nuhag-php/bibtex/
open\_files/fe03-1\_modspa03.pdf}.

\bibitem{GJ} J. Ginibre, G. Velo, Generalized Strichartz
inequalities for the wave equations. J.Funct. Anal., {\bf 133}
(1995), 50--68.

\bibitem{GV2} J. Ginibre, G. Velo, The Cauchy problem in local spaces for the complex
Ginzburg--Landau equation II, Contraction Methods,
 Commun. Math. Phys.,  {\bf 187} (1997),  45--79.

 \bibitem{GL} V. Ginzburg and L. Landau,
 On the theory of superconductivity,
 Zh. Eksp. Fiz.,  {\bf 20} (1950), 1064; English transl. in: Men of Physics:
 L. D. Landau, Vol. L, D. Ter Haar(Ed.), New York: Pergamon Press,
 1965, 546--568.


\bibitem{GW1}   B. L. Guo and B. X. Wang, Finite dimensional behavior for the derivative Ginzburg-Landau
equation in two spatial dimensions, Physica D.,  {\bf 89}(1995),
83--99.

\bibitem{12}     H. Gao and J. Duan, On the initial value problem for the generalized
2D Ginzburg-Landau equation, J. Math. Anal. Appl., {\bf 216}
(1997), 536--548.


\bibitem{FDF} L. Han, B. Wang,
Global wellposedness and limit behavior for the generalized
finite-depth-fluid equation with small critical data, J.
Differential Equation., {\bf 245} (2008), 2103--2144.

\bibitem{Huo} Z. Huo, Y. Jia, Well-posedness and inviscid limit behavior of solution for the
generalized 1D Ginzburg--Landau equation. J. Math. Pures Appl.,
 {\bf 92}(2009) 18--51.

\bibitem{HH} C. Huang, B. Wang, Inviscid limit for the
energy-critical complex Ginzburg-Landau equation, Journal of
Functional Analysis., {\bf 255} (2008), 681--725.

\bibitem{NH}  N. Hayashi,  The initial value problem for the derivative
nonlinear Schr\"odinger equation in the energy space. Nonlinear
Anal., TMA., {\bf 20} (1993)823--833.

\bibitem{NH2}N. Hayashi, T. Ozawa, Remarks on nonlinear
Schr\"odinger equations in one space dimension. Differential and
Integral Equations.,  {\bf 2}(1994),   453--461.

\bibitem{ICE} A. Ionsescu and C. E. Kenig, Low-regularity
Schr\"odingger maps, II:Global well posedness in dimension $d\geq
3$, Commun. Math. Physics.,   {\bf 271}  (2007), 523--559.

\bibitem{KT}  T. Kato, Nonlinear Schr\"odinger Equation . II.
$H^s$-solutions and unconditional well-possedness. J.Anal. Math.,
{\bf 67 }(1995), 281--306.

\bibitem{KTao} M. Keel, T. Tao, Endpoint Strichartz estimates. Amer.
J. Math., {\bf  120} (1998),  955--980.

\bibitem{KPVL} C. E. Kenig, G. Ponce, L. Vega, Oscillatory integrals and
regularity of dispersive equations, Indiana Univ. Math. J.,  {\bf
40} (1991), 32--69.

\bibitem{KPVH} C. E. Kenig, G. Ponce, L. Vega,  Higher-order
nonlinear dispersive equations. Proc. Amer. Math. Soc., {\bf 122}
(1994),  157--166.

\bibitem{14} C. E. Kenig, G. Ponce, L. Vega, Small solutions to nonlinear
Schrodinger equation, Ann. Inst. H. Poincar\'{e}. Sect. C., {\bf 10}
(1993) 255--288.

\bibitem{15} C. E. Kenig, G. Ponce, L. Vega, Smoothing effects and local
existence theory for the generalized nonlinear Schr\"odinger
equations, Invent. Math., {\bf 134} (1998) 489--545.



\bibitem{KT} Kato, T. Nonlinear Schr\"odingger equation . II.
$H^s$-solutions and unconditional well-possedness. J.Anal. Math.,
{\bf 67} (1995), 281--306.

\bibitem{18} S. Klainerman, Long-time behavior of solutions to nonlinear
evolution equations, Arch. Ration. Mech. Anal., {\bf 78} (1982)
73--98.

\bibitem{19} S. Klainerman, G. Ponce, Global small amplitude solutions to
nonlinear evolution equations, Comm. Pure Appl. Math., {\bf 36}
(1983) 133--141.

\bibitem{LG}  Y. Li, B. Guo, Global existence of solutions to the 2D
Ginzburg¨CLandau equation, J. Math. Anal. Appl., {\bf 249} (2000)
412--432.





\bibitem{MN} S. Machihara, Y. Nakamura, The inviscid limit for the complex
Ginzburg-Landau equation, JMAA.,  {\bf 281} (2003), 552--564.

\bibitem{22}  T. Ozawa, J. Zhai, Global existence of small
classical solutions to nonlinear Schr\"odinger equations, Ann. Inst.
H. Poincar\'{e} Anal. Non Linear., {\bf 25} (2008) 303--311.

\bibitem{Li-Po93} F. Linares and G. Ponce, On the Davey-Stewartson systems,   Ann. Inst.
H. Poincar\'e,  Anal. NonLin\'eaire.,  {\bf 10} (1993), 523--548.


\bibitem{stein1} E. M. Stein,. Harmonic Analysis. Princeton
University Press, 1993.

\bibitem{stein2} E. M. Stein and Guido Weiss,. Introduction to
Fourier Analysis on Euclidean Spaces. Princeton University Press,
1971.

\bibitem{ATA} A. Stefanov, On quadratic derivative Schr\"odingger
equations in one space dimension, Transactions of American
Mathematical Sociaty., {\bf 359} (2007), 3589--3607.

\bibitem{WH}  B. Wang, L. Han, C. Huang, Global
Well-posedness and Scattering for the Derivative Nonlinear
Schr\"odingger Equation with Small Rough Data,  Ann. Inst. H.
Poincar\'{e}, Anal. Nonlin\'{e}aire.26(2009),2253-2281.

\bibitem{Wang} B. Wang, The limit behavior of solutions for the Cauchy
problem of the Complex Ginzburg-Landau equation, Commu. Pure. Appl.
Math., {\bf 55}  (2002), 0481--0508.

\bibitem{WGZ} B. Wang, B. Guo, L. Zhao, The global well-posedness and spatial decay of solutions for the
derivative complex Ginzburg--Landau equation in $H^1$. Nonlinear
Anal., {\bf 57} (2004) 1059--1076.

\bibitem{WaHe}  B. X.  Wang and H. Hudzik, {\rm The global Cauchy problem for the
NLS and NLKG with small rough data,} J. Differential Equations.,
{\bf 231} (2007), 36--73.

\bibitem{Wa2} B. X. Wang, Y. D. Wang, The inviscid limit for the derivative
Ginzburg-Landau equations, J.  Math. Pures Appl., {\bf 83} (2004),
477-502.

\bibitem{WU} J. H. Wu,   The inviscid limit of the complex
Ginzburg-Landau equation, J. Differential Equations., {\bf 142}
(1998), 413--433.
\end{thebibliography}
\end{document}